\documentclass[11pt]{article}
\usepackage{latexsym,amssymb}
\usepackage{amsmath}
 \usepackage{paralist}

\newtheorem{Theorem}{\bf Theorem}[section]
\newtheorem{Lemma}{\bf Lemma}[section]
\newtheorem{Proposition}{\bf Proposition}[section]
\newtheorem{Corollary}{\bf Corollary}[section]
\newtheorem{Remark}{\bf Remark}[section]
\newtheorem{Example}{\bf Example}[section]
\newtheorem{Definition}{\bf Definition}[section]

\newenvironment{lemma}{\begin{Lemma}$\!\!\!$}{\end{Lemma}}

\newenvironment{definition}{\begin{Definition}$\!\!\!$}{\end{Definition}}

\numberwithin{equation}{section}

\begin{document}

\title{Higher order asymptotic expansions to the solutions for\\
 a nonlinear damped wave equation
\quad}

\author{Tatsuki Kawakami\\
Department of Mathematical Sciences, Osaka Prefecture University,\\
Sakai 599-8531, Japan\\
(e-mail address : kawakami@ms.osakafu-u.ac.jp)
\quad\vspace{8pt}\\
and\vspace{8pt}\\
Hiroshi Takeda\\
Faculty of Engineering, Fukuoka Institute of Technology, \\
Fukuoka, 811-0295 Japan \\ 
(e-mail address : h-takeda@fit.ac.jp)
}
\date{}
\maketitle

\begin{abstract}
We study the Cauchy problem for a nonlinear damped wave equation.
Under suitable assumptions for the nonlinearity and the initial data,
we obtain the global solution which satisfies weighted $L^1$ and $L^\infty$ estimates.
Furthermore, we establish the higher order asymptotic expansion of the solution.
This means that we construct the nonlinear approximation of the global solution 
with respect to the weight of the data.
 Our proof is based on the approximation formula of the linear solution, which is given in \cite{T01},
 and the nonlinear approximation theory for a nonlinear parabolic equation developed by \cite{IKK}.
\end{abstract}
\section{Introduction}

We consider the Cauchy problem 
for a nonlinear damped wave equation,
\begin{equation}
\label{eq:1.1}
\left\{
\begin{array}{l}
\displaystyle
\partial_t^2u-\Delta u+\partial_tu=F(u)\quad \mbox{in}\quad {\mathbb R}^N\times(0,\infty),
\vspace{5pt}\\
\displaystyle
u(x,0)=u_0(x),\quad\partial_tu(x,0)=u_1(x)\quad \mbox{in}\quad{\mathbb R}^N,
\end{array}
\right.
\end{equation}
where $N = 1,2,3$, and $\partial_t=\partial/\partial t$.
We assume that the nonlinear term $F\in C({\mathbb R})$ satisfies 
\begin{align}
\label{eq:1.2}
|F(\xi)|&\le C|\xi|^p,\qquad \xi\in{\mathbb R},
\\
\label{eq:1.3}
|F(\xi)-F(\eta)|&\le C(|\xi|^{p-1}+|\eta|^{p-1})|\xi-\eta|,\qquad \xi,\eta\in{\mathbb R},
\end{align}
for some constants $p>1$ and $C>0$, which are independent of $\xi$ and $\eta$.
The typical examples of our nonlinear terms are given by 
\begin{align} \label{eq:1.4}
F(\xi)=\pm|\xi|^{p-1} \xi\quad\mbox{or}\quad \pm |\xi|^{p}\qquad (p>1).
\end{align}
The aim of this paper is to study the large time behavior of the solution to \eqref{eq:1.1}.
More precisely, we show the nonlinear approximation of the solution to \eqref{eq:1.1} 
with respect to the order of the moment for the initial data.
This point of view is shared with Ishige and  the first author of the paper \cite{IK} 
studying a nonlinear heat equation.

In the classical paper \cite{M}, 
Matsumura considered the Cauchy problem of nonlinear wave equations with dissipation terms.
His main tools in the proof are the estimates for the solutions to a linear damped wave equation
\begin{equation}
\label{eq:1.5}
\left\{
\begin{array}{l}
\displaystyle
\partial_t^2v-\Delta v+\partial_tv=0\quad \mbox{in}\quad {\mathbb R}^N\times(0,\infty),
\vspace{5pt}\\
\displaystyle
v(x,0)=v_0(x),\quad\partial_tv(x,0)=v_1(x)\quad \mbox{in}\quad{\mathbb R}^N.
\end{array}
\right.
\end{equation}
Especially, he prove the decay estimates for the solution to \eqref{eq:1.5} by the Fourier splitting method.
Beginning this paper, 
many authors showed the large time behavior of the solution to \eqref{eq:1.5}
(see, e.g., \cite{BF, Li, YM}).
In \cite{OZP},
Orive, Zuazua and Pozato obtained 
the higher order asymptotic expansion of the solution to \eqref{eq:1.5} 
for the variable coefficient setting in $L^{2}$ base framework.
Furthermore, 
the decomposition of the solution into the solutions of heat equations 
and wave equations are proposed
(see, e.g., \cite{MN} for $N=1$, \cite{HO} for $N=2$,
\cite{N} for $N=3$ and \cite{Na} for $N \ge 4$).

On the other hand,
the Cauchy problem \eqref{eq:1.1} with \eqref{eq:1.4} has been studied by many mathematicians
from various points of view.
Especially, for the focusing case $F(\xi)=|\xi|^p$,
it is well-known that the growth order $p=1+2/N$ is critical situation 
for the existence of the global solution to \eqref{eq:1.1}.
In \cite{TY},
Todorova and Yordanov proved that,
if $p>1+2/N$, then there exists a unique global solution for the small compactly supported data, 
and if $1< p < 1+2/N$, then the solution blows up in a finite time for small data 
(see also \cite{LZ} for blow up results).
Zhang \cite{Z} showed the small data blow-up results including $p=1+2/N$.
Ikehata and Tanizawa \cite{IT} obtained the global existence results for $p>1+2/N$ 
under the non-compact supported data assumption
(see,e.g., \cite{IMN} for the case $F(\xi)=|\xi|^{p-1}\xi$).
For the defocusing case $F(\xi) = -|\xi|^{p-1} \xi$ or  $F(\xi) = -|\xi|^{p}$,
it is also well-investigated,
and it is well known that there exists a unique global solution with decay property for all $p>1$ 
if data has sufficient regularity
(see, e.g., \cite{HKN2, INZ, KNO}).

For the asymptotic profiles of the solution to damped wave equations, 
so-called diffusion phenomena is shown by many authors.
Among others,
Gallay and Raugel \cite{GR} proved that global solutions of nonlinear damped wave equation 
behaves like those of nonlinear heat equations with suitable data,
including more general nonlinearity for $N=1$.
In \cite{K},
Karch proved the approximation of the solution to \eqref{eq:1.1} by the heat kernel for $p \ge 1+ 4/N$.
After that,
Nishihara \cite{N} proved it  for $p>1+2/N$ when $N=3$.
(See also \cite{MN} for $N=1$, \cite{HO} for $N=2$, \cite{Na} for $N=4,5$ 
and \cite{HKN1} for $N \ge 1$).
In \cite{KU},
Ueda and the first author of this paper
obtained the second order nonlinear approximation of the solution to \eqref{eq:1.1} for $p>1+2/N$.
Recently, 
the second author of this paper \cite{T02} proved 
the $K+1$ th order expansion of the solution by the series of the heat kernel when $p> 1+ (K+2)/N$.

We should also state the topic on the recent progress of 
the diffusion phenomena in damped wave equations.
One of our motivation here, the precise description of the large time behavior of the solution, 
is shared with the related results for the following equations:
\begin{itemize}
\item The variable coefficient damping case (see, e.g., \cite{LNZ1, LNZ2, N1, N3, Wa, W1, W2, Ya})
$$
\partial_t^2v-\Delta v+ b(t,x) \partial_tv=0 \quad \mbox{in}\quad {\mathbb R}^N\times(0,\infty);
$$
\item The structural damed wave equations
(see, e.g., \cite{I1, ITY, JS, P, S} for $\sigma=1$ and \cite{DE1, DE2, DR} for $0<\sigma\le1$)
$$
\partial_t^2v-\Delta v+(-\Delta)^{\sigma} \partial_tv=0\quad \mbox{in}\quad {\mathbb R}^N\times(0,\infty).
$$
\end{itemize}
We remark that based on the linear estimates, 
the authors also treated the nonlinear perturbation 
by the various methods to have the diffusion phenomena.  
\vspace{5pt}

As we seen the above, 
there are many results for the diffusion phenomena of the dissipative type wave equations. 
On the other hand, there are few studies concerning the higher order asymptotic expansion 
of the global solutions of \eqref{eq:1.1}
except for the results \cite{GR} (1-d case), \cite{KU} (up to second order expansion) 
and \cite{T02} (under the regularity for the nonlinearity).
Roughly speaking, the difficulty stems from
the construction of the approximation functions as $t \to \infty$ 
and the singularity from the high frequency part.
To avoid the difficulty, 
the previous results need to restrict the order of the expansion or to only treat the smooth nonlinearity. 
In this paper, we established the asymptotic expansion of the solution to \eqref{eq:1.1} 
for \eqref{eq:1.4} with $p>1+2/N$ $(N=1,2)$,  
$p\ge 2$ $(N=3)$ up to the suitable order,
which depends on the order of the moment for the initial data.
Our proof is based on the higher order asymptotic expansion formula of the linear solution 
by the solution of the heat equation,
which is  shown by the second author of this paper in \cite{T01}, 
and nonlinear approximation technique for nonlinear parabolic equations 
developed by Ishige, Kobayashi and the first author of the paper in \cite{IKK}.      
Our new ingredient here is to show the weighted $L^{\infty}$ estimates 
for the solution to \eqref{eq:1.5},
which are important to apply the iteration scheme proposed in \cite{IKK} to our problem \eqref{eq:1.1}.
Furthermore, 
our results imply not only the detailed profile of the solution to \eqref{eq:1.1} for large $t$, 
but also sharp decay estimates in each expansion order.
\vspace{5pt}

This paper is organized as follows. 
In section 2, we prepare notations which used throughout this paper
and we state the main results in this paper.
Furthermore we mention important remarks on main results.
Section 3 presents some preliminaries.
Section 4 is devoted to the study of the weighted $L^1$ and $L^{\infty}$ estimates 
for the solutions to \eqref{eq:1.5}.
In section 5 and section 6, our main results are proved. 

\section{Main Results}

\subsection{Notation} 
To state our results precisely, 
we summarize notion and notation.
For $k\ge0$, we denote $[k]$ the integer satisfying $k-1<[k]\le k$.
Let ${\mathbb N}_0={\mathbb N}\,\cup\,\{0\}$, ${\mathbb M}={\mathbb N}_0^N$,
and $G$ be the $N$-dimensional heat kernel, that is,
\begin{equation}
\label{eq:2.1}
G(x,t) :=(4\pi t)^{-N/2}\exp\left(-\frac{|x|^2}{4t}\right).
\end{equation}
Furthermore,
for any $\phi\in L^q({\mathbb R}^N)$ with $q\in[1,\infty]$, 
we denote by $e^{t\Delta}\phi$ 
the unique bounded solution for the Cauchy problem of the heat equation with the initial datum $\phi$, 
that is, 
\begin{equation}
\label{eq:2.2}
(e^{t\Delta}\phi)(x):=\int_{{\mathbb R}^N}G(x-y,t)\phi(y)dy.
\end{equation}
For any $\alpha=(\alpha_1,\cdots,\alpha_N)\in{\mathbb M}$, 
we put 
\begin{align}
&
|\alpha|:=\displaystyle{\sum_{i=1}^N}|\alpha_i|,\quad
\alpha!:=\prod_{i=1}^N\alpha_i!,\quad
x^\alpha:=\prod_{i=1}^N x_i^{\alpha_i},\quad
\partial_x^\alpha:=
\frac{\partial^{|\alpha|}}{\partial x_1^{\alpha_1}\cdots\partial x_N^{\alpha_N}},
\vspace{5pt}
\notag\\
&
g(x,t):=G(x,1+t),
\qquad
g_\alpha(x,t):=\frac{(-1)^{|\alpha|}}{\alpha!}\partial_x^\alpha g(x,t).
\notag
\end{align}
Let ${\mathbb M}_k:=\left\{\alpha\in{\mathbb M}:|\alpha|\le k\right\}$ for $k\ge 0$. 
For any $\alpha=(\alpha_1,\dots,\alpha_N)$, $\beta=(\beta_1,\dots,\beta_N)\in{\mathbb M}$, 
we say 
$$
\alpha\le\beta
$$
if $\alpha_i\le\beta_i$ for all $i\in\{1,\dots,N\}$.  
For $\ell$, $m$ and $n \in \mathbb{N}_0$,
we put
$$
\phi_{\ell}(r)
 := 
\left(
\frac{2}
{1 + \sqrt{1 -4r}}
\right)^{2\ell}, \qquad 
\psi(r) 
:= \frac{1}{\sqrt{1 - 4r}},
$$
and
\begin{equation} 
\label{eq:2.3}
\Phi_{\ell,m} := \frac{1}{\ell! m!} \frac{d^{m}}{dr^{m}} \phi_{\ell}(r) \biggr|_{r=0}, \qquad
\Psi_{n} := \frac{1}{\ell !} \frac{d^{n}}{dr^{n}} \psi(r) \biggr|_{r=0}.
\end{equation}
For any $q\in[1,\infty]$,
we denote by $L^q$ and $\|\cdot\|_{L^q}$ the usual $L^q({\mathbb R}^N)$ space and its norm,
respectively.
Let $\ell\in {\mathbb N}_0$.
Then $W^{\ell,q}$ denotes the Sobolev space of $L^q$ functions,
equipped with the norm
$$
\|\phi\|_{W^{\ell,q}}:=\left(\sum_{|\alpha|\le \ell}\|\partial_x^\alpha \phi\|_{L^q}^q\right)^{1/q}.
$$
For any $k \ge 0$ and $q\in\{1,\infty\}$, we denote by $L^q_k$ the functional space  
$L^q({\mathbb R}^N,(1+|x|)^k dx)$, and put 
$$
\|\phi\|_{L^q_k}=\|\phi\|_{L^q((1+|x|)^k dx)},
\qquad
\|\phi\|_{W^{\ell,q}_k}=\sum_{|\alpha|\le \ell}\|\partial_x^\alpha \phi\|_{L^q_k}.
$$
Here we often identify $W^{0,1}_{k}=L^{1}_{k}$.
Throughout the present paper,
$C$ denotes a various generic positive constant.
\vspace{5pt}

Let us give the definition of the solution for the Cauchy problem \eqref{eq:1.1}.
\begin{definition}
\label{Definition:2.1}
Let $u\in C([0,\infty):L^1) \cap L^\infty(0,\infty:L^\infty)$
and $F\in C({\mathbb R})$.
Then the function $u$ is said to be a solution for the Cauchy problem \eqref{eq:1.1} 
if there holds 
\begin{equation}
\label{eq:2.4}
u(x,t)=(K_0(t)u_0)(x)+\bigg(K_1(t)\left(\frac{1}{2}u_0+u_1\right)\bigg)(x)+\int_0^t(K_1(t-s)F(u(s)))(x)ds
\end{equation}
for all $(x,t)\in{\mathbb R}^N\times(0,\infty)$,
where the evolution operators $K_0(t)$ and $K_1(t)$ of the linear damped wave equation
are given by
\begin{align} 
    (K_{0}(t)\phi)(x) & := \mathcal{F}^{-1} \left[e^{-\frac{t}{2}}
	\cos\left(t \sqrt{|\xi|^{2}- \frac{1}{4}}\right)
	\mathcal{F}[\phi] \right](x),
	\notag\\
\label{eq:2.5}
	(K_{1}(t)\phi)(x) & := \mathcal{F}^{-1}
    \left[e^{-\frac{t}{2}}
	\frac{\sin
    \left(t \sqrt{|\xi|^{2}- \frac{1}{4}}\right)}
        {\sqrt{|\xi|^{2}- \frac{1}{4}}}
	 \mathcal{F}[\phi] \right](x). 
\end{align}
Here we denote the Fourier and Fourier inverse transform by
$\mathcal{F}$ and  $\mathcal{F}^{-1}$, respectively.
\end{definition}

Let $M_\alpha(f,t)$ be the constant defined inductively (in $\alpha$) by
\begin{equation}
\label{eq:2.6}
\begin{array}{l}
\displaystyle{M_0(f,t):=\int_{{\mathbb R}^N}f(x)\,dx},\quad\mbox{if}\quad\alpha=0,\vspace{7pt}\\
\displaystyle{M_\alpha(f,t):=\int_{{\mathbb R}^N} 
x^\alpha f(x)\,dx-\sum_{\beta\le\alpha,\beta\not=\alpha}M_\beta(f,t) 
\int_{{\mathbb R}^N}x^\alpha g_\beta(x,t)\,dx}\quad\mbox{if}\quad\alpha\not=0.
\end{array}
\end{equation}
By \eqref{eq:2.3} and \eqref{eq:2.6}
we introduce the function $U_{lin}(t)=U_{lin}(x,t)$ by 
\begin{equation}
\label{eq:2.7}
\begin{split}
&U_{\rm lin}(x,t)
\\
&=
 \sum_{\ell=0}^{[K/2]} 
 \sum_{m=0}^{[K/2]-\ell} \Phi_{\ell, m}
 \Bigg\{\frac{1}{2}
 \sum_{|\alpha| \le K-2(\ell+m)}
 M_{\alpha}(u_0,0)
(-t)^{\ell} 
(-\Delta)^{2\ell+m}g_\alpha(x,t)
\\
&\qquad +
\sum_{n=0}^{[K/2]-\ell-m}
\sum_{|\alpha| \le K-2(\ell+m +n)} 
\Psi_{n} 
\left(\frac{1}{2}M_{\alpha}(u_0,0) 
+M_{\alpha}(u_1,0) \right)\times
\\
&
\qquad\qquad\qquad\qquad\qquad\qquad\qquad\qquad\qquad
\times
(-t)^{\ell} (-\Delta)^{2\ell+m+n}g_\alpha(x,t)
\Bigg\}.
\end{split}
\end{equation}
Following the notation of \cite{IIK},
we also introduce the linear operator $P_{K}$ on $L_{K}^{1}$ by 
$$
[P_K(t)f](x):=f(x)-\sum_{|\alpha|\le K}M_\alpha(f,t) g_\alpha(x,t),
$$
for $K\ge 0$ and $t>0$ and $f\in L^1_K$.
Then the operator $P_K(t)$ has the following property,
\begin{equation}
\label{eq:2.8}
\int_{{\mathbb R}^N}x^\alpha[P_k(t)f](x)\, dx=0,
\qquad t>0.
\end{equation}
for any $\alpha\in {\mathbb M}_K$.
This property plays an important role of deriving our main results
(for the detail, see also {\rm\cite{IIK, IK, KU}}). 

With the notation and the notion above, we can define the sequence of functions $U_j=U_j(x,t)$ 
inductively:
Let $K\ge0$ and $(u_{0}, u_{1}) \in W^{[N/2],1}_{K} \cap W^{[N/2], \infty}\times 
L^{1}_{K} \cap L^{\infty}$.
For any $j \in {\mathbb N}$, 
\begin{equation}
\label{eq:2.9}
\begin{split}
&
U_0(t):=U_{\rm lin}(t)+\sum_{|\alpha|\le K}\int_0^tK_1(t-s)M_\alpha(F(u(s)),s)g_\alpha(s)\,ds,
\\
&
U_j(t):=
U_0(t)+\int_0^tK_1(t-s)P_K(s)F_{j-1}(s)\,ds,
\end{split}
\end{equation}
where $F_j(x,t):=F(U_j(x,t))$.
Here $K_1(t)$ and $U_{\rm lin}$ are given in \eqref{eq:2.5} and \eqref{eq:2.7}, respectively.

\subsection{Main Theorems}
Now we are ready to treat our main results.
Our first result gives the sufficient condition for the existence of global solution 
satisfying the weighted $L^1$ and $L^{\infty}$ estimates with suitable decay property. 
\begin{Theorem}
\label{Theorem:2.1}
Let $K\ge0$ and $(u_{0}, u_{1}) \in W^{[N/2],1}_{K} \cap W^{[N/2], \infty}_{K}\times 
L^{1}_{K} \cap L^{\infty}_{K}$.
Assume that $F\in C({\mathbb R})$ satisfies \eqref{eq:1.2} and \eqref{eq:1.3} with 
\begin{equation}
\label{eq:2.10}
p>p_*:=1+\frac{2}{N}\quad (N=1,2),\qquad
p\ge2\quad (N=3).
\end{equation}
Then there exists a positive constant $\varepsilon$ such that
if $E_0:=\|u_0\|_{W^{[N/2],\infty}}+\|u_0\|_{W^{[N/2],1}}
+\|u_1\|_{L^\infty}+\|u_1\|_{L^1} \le\varepsilon$, then the Cauchy problem \eqref{eq:1.1} admits 
a unique global solution $u$ of \eqref{eq:1.1} in the class
$$C([0,\infty): L^1)\cap L^\infty(0,\infty:L^\infty)$$
satisfying
\begin{equation}
\label{eq:2.11}
\|u(t)\|_{L^q}\le CE_0(1+t)^{-\frac{N}{2}(1-\frac{1}{q})}, \quad t\ge0,
\end{equation}
for any $q\in[1,\infty]$.
Moreover, $$u \in C([0,\infty): L^1_{K})\cap L^\infty(0,\infty:L^\infty_{K})$$
and 
\begin{align}
\label{eq:2.12}
&
\|u(t)\|_{L^1_k}\le CE_K(1+t)^{\frac{k}{2}},\qquad t\ge0,\\
&
\|u(t)\|_{L^\infty_k}\le CE_K(1+t)^{\frac{k-N}{2}},\qquad t\ge0,
\label{eq:2.13}
\end{align}
for any $k\in[0,K]$,
where 
\begin{equation}
\label{eq:2.14}
E_K:=\max_{k_i\in[0,K],i\in\{1,2,3,4\}}\Bigg\{\|u_0\|_{W^{[N/2],\infty}_{k_1}}+\|u_0\|_{W^{[N/2],1}_{k_2}}
+\|u_1\|_{L^\infty_{k_3}}+\|u_1\|_{L^1_{k_4}}\Bigg\}.
\end{equation}
\end{Theorem}
\begin{Remark}
Here we note that Theorem~$\ref{Theorem:2.1}$ states that to obtain the suitable estimates 
\eqref{eq:2.12} and \eqref{eq:2.13}, we only assume that the smallness of $E_{0}$, not $E_{K}$.
\end{Remark}

Our second result is the nonlinear approximation of the solutions to \eqref{eq:1.1}.
In other words, 
the functions $U_{j}$ defined by \eqref{eq:2.9} are nonlinear approximation 
of the global solution to \eqref{eq:1.1}.
(See also Remark~\ref{Remark:3.2}.)

\begin{Theorem}
\label{Theorem:2.2}
Let $j\in{\mathbb N}_0$, and let $U_j$ be the functions given in \eqref{eq:2.9}.
Put
\begin{equation}
\label{eq:2.15}
A:=\frac{N}{2}(p-1)-1>0.
\end{equation}
Assume that there exists a unique global solution to \eqref{eq:1.1} satisfying 
\eqref{eq:2.11}, \eqref{eq:2.12} and \eqref{eq:2.13}.
Then, for any $q\in[1,\infty]$, $k\in[0,K]$ and $\gamma\in\{1,\infty\}$, 
\begin{equation}
\label{eq:2.16}
\sup_{t>0}\,(1+t)^{\frac{N}{2}(1-\frac{1}{q})}\|U_j(t)\|_{L^q}
+\sup_{t>0}\,(1+t)^{-\frac{k}{2}+\frac{N}{2}(1-\frac{1}{\gamma})}\|U_j(t)\|_{L^\gamma_k}
<\infty
\end{equation}
and 
\begin{equation}
\label{eq:2.17}
 t^{\frac{N}{2}(1-\frac{1}{q})}\left\|u(t)-U_j(t)\right\|_{L^q}
 =\left\{
\begin{array}{ll}
o(t^{-\frac{K}{2}})+O(t^{-(j+1)A})  & \mbox{if}\quad (j+1)A\not=K/2,\vspace{3pt}\\
O(t^{-\frac{K}{2}}\log t)  & \mbox{if}\quad (j+1)A=K/2,\\
\end{array}
\right.
\end{equation}
as $t\to\infty$. 
\end{Theorem}
\begin{Remark}
Theorem~$\ref{Theorem:2.2}$ shows that once we have a unique solution to \eqref{eq:1.1} satisfying
the estimates \eqref{eq:2.11}, \eqref{eq:2.12} and \eqref{eq:2.13}, without smallness of the data, 
then we see that the sequence of the functions $U_{j}$ is well-defined 
\eqref{eq:2.16} suitable sense and 
we have the asymptotic behavior of the solutions as $t \to \infty$ \eqref{eq:2.17}.
\end{Remark}
%

\section{Preliminaries}

\subsection{Solutions of the heat equation}
In this subsection we recall some preliminary results
on the behavior of solutions for the heat equation and the operator $P_k(t)$. 
\vspace{5pt}

Let $\alpha\in{\mathbb M}$ and $G$ 
be the function given in \eqref{eq:2.1}.
Then 
we have
\begin{equation}
\label{eq:3.1}
|\partial^\alpha_xG(x,t)|\le Ct^{-\frac{N+|\alpha|}{2}}
\left[1+\left(\frac{|x|}{t^{1/2}}\right)^{|\alpha|}\right]
\exp\left(-\frac{|x|^2}{4t}\right)
\end{equation} 
for all $(x,t)\in{\mathbb R}^N\times(0,\infty)$. 
This inequality yields the inequalities
\begin{align}
&
\label{eq:3.2}
\|g_\alpha(t)\|_{L^q}\le C (1+t)^{-\frac{N}{2}(1-\frac{1}{q})-\frac{|\alpha|}{2}},
\quad t>0,
\\
&
\label{eq:3.3}
\|g_\alpha(t)\|_{L^\gamma_k}\le C (1+t)^{\frac{k-|\alpha|}{2}-\frac{N}{2}(1-\frac{1}{\gamma})},
\quad t>0,
\end{align} 
for any $q\in[1,\infty]$, $k\ge 0$ and $\gamma\in\{1,\infty\}$.
Furthermore, applying the Young inequality to \eqref{eq:2.2}  with the aid of \eqref{eq:3.1},
for any $\alpha\in{\mathbb M}$ and $1\le r\le q\le\infty$, we have
$$
\|\partial_x^\alpha e^{t\Delta}\phi\|_{L^r}
\le C t^{-\frac{N}{2}(\frac{1}{r}-\frac{1}{q})-\frac{|\alpha|}{2}}\|\phi\|_{L^r},\quad t>0.
$$
In particular, we obtain
\begin{equation*}
\|e^{t\Delta}\phi\|_{L^r}\le\|\phi\|_{L^r},\quad t>0.
\end{equation*}

We recall the following lemma,
which is useful in our study for the asymptotic expansion of solutions. 

%
\begin{Lemma}
\label{Lemma:3.1}
{\rm(\cite[Proposition 3.1]{IKK}.)}
Let $K\ge0$. 
Then the following holds.
\vspace{3pt}
\newline
{\rm (i)}
For any $k\in[0,K]$, 
$$
\int_{{\mathbb R}^N}|x|^k|e^{t\Delta}\phi(x)|\,dx
\le Ct^{-\frac{K-k}{2}}\int_{{\mathbb R}^N}|x|^K|\phi(x)|\,dx,
\quad t>0,
$$
for all $\phi\in L^1_K$ satisfying 
\begin{equation}
\label{eq:3.4}
\int_{{\mathbb R}^N}x^\alpha\phi(x)\,dx=0,
\quad \alpha\in{\mathbb M}_K.
\end{equation}
{\rm (ii)} 
For any $k\in[0,K]$,  
$$
\lim_{t\to\infty}t^{\frac{K-k}{2}}\int_{{\mathbb R}^N}|x|^k|e^{t\Delta}\phi(x)|\,dx=0
$$
for $\phi\in L^1_K$ satisfying \eqref{eq:3.4}.
\end{Lemma}

Let us mention the important fact to show the proof of Theorem \ref{Theorem:2.2}.
The point of the following Proposition \ref{Proposition:3.1} is in the assertion that 
$e^{t \Delta} \phi$ is well-approximated by the sequence of $g_{\alpha}(t)$, 
not $\partial_{x}^{\alpha} G(x,t)$.  
\begin{Proposition}
{\rm(\cite[Theoreme~1.1]{IKK})}
\label{Proposition:3.1}
Let $K\ge0$ and $\phi\in L^1_K$.
Then, for any $j\in{\mathbb N}_0$ and $q\in[1,\infty]$,
$$
t^{\frac{N}{2}(1-\frac{1}{q})+\frac{j}{2}}
\left\|\nabla^j\left[e^{t\Delta}\phi-\sum_{|\alpha|\le K}M_\alpha(\phi,0)g_\alpha(t)\right]\right\|_{L^q} =
o(t^{-\frac{K}{2}})
$$
as $t \to \infty$.
\end{Proposition}

At the end of this subsection we give a lemma on the estimate of the function $P_K(t)f(t)$. 

\begin{lemma}
\label{Lemma:3.2}
Let $K\ge0$ and $1\le q\le\infty$.
Let $f$ be a measurable function in ${\mathbb R}^N\times(0,\infty)$ such that 
\begin{equation}
\label{eq:3.5}
\begin{split}
E_{K,q}[f](t):=
&
(1+t)^{\frac{K}{2}}\left[(1+t)^{\frac{N}{2}(1-\frac{1}{q})}\|f(t)\|_{L^q}
+\|f(t)\|_{L^1}\right]
\\
&
\qquad
+\|f(t)\|_{L^1_K}+(1+t)^{\frac{N}{2}}\|f(t)\|_{L^\infty_K}
\in L^\infty(0,T)
\end{split}
\end{equation}
for any $T>0$. 
Then, for any $1\le r\le q$, $0\le k\le K$ and $\gamma\in\{1,\infty\}$,
\begin{equation}
 |M_\alpha(f(t),t)|\le 
 C(1+t)^{-\frac{K-|\alpha|}{2}}E_{K,q}[f](t),\quad \alpha\in{\mathbb M}_K,
  \label{eq:3.6}
\end{equation}
and
\begin{equation}
\begin{split}
&
(1+t)^{\frac{N}{2}(1-\frac{1}{q})}\|P_K(t)f(t)\|_{L^q}
+(1+t)^{-\frac{k}{2}+\frac{N}{2}(1-\frac{1}{\gamma})}\|P_K(t)f(t)\|_{L^\gamma_k}
\\
&
\qquad\qquad\qquad\qquad\qquad\qquad\qquad\qquad\qquad
\le C(1+t)^{-\frac{K}{2}}E_{K,q}[f](t)
\end{split}
\label{eq:3.7}
\end{equation}
for almost all $t>0$. 
\end{lemma}
\noindent
{\bf Proof.}
Applying the same arguments as in the proof of \cite[Lemma~2.2]{IKK},
we can easily prove this lemma.
So we omit the proof.
$\Box$
\vspace{8pt}
%
\subsection{Solutions of the damped wave equation}

In this subsection, we recall some preliminary results on 
the properties of the fundamental solutions 
for the linearized damped equation \eqref{eq:1.5}.
We first give the well-known representation formulas of 
evolution operators $K_{1}(t)\phi$ and $K_{0}(t)\phi$ 
(see, e.g., \cite{CH, KU, N, T01}).
For $N=1,2,3$, we have
\begin{equation*}
  (K_1(t)\phi)(x) =
\begin{cases}
& \dfrac{e^{-\frac{t}{2}}}{2}
           \displaystyle\int_{|z|\le t}
                I_0 \left({\sqrt{t^2-|z|^2}}/2 \right) \,\phi(x+z)\,dz
\quad (N=1),
\vspace{5pt}\\ 
  & \dfrac{e^{-\frac{t}{2}}}{2 \pi}
           \displaystyle\int_{|z|\le t}
                \dfrac{\cosh \left( {\sqrt{t^2-|z|^2}}/2 \right)}
                     {\sqrt{t^2-|z|^2}}\,
                \phi(x+z)\,dz \quad (N=2),
                \vspace{5pt}\\
& \dfrac{e^{-\frac{t}{2}}}{4 \pi t}\,
           \partial_t
           \displaystyle\int_{|z|\le t}
                I_0 \left( {\sqrt{t^2-|z|^2}}/2 \right)\,
                \phi(x+z)\,dz \quad (N=3), 
\end{cases}
\end{equation*}
where $I_{\nu}$ is the modified Bessel function of order 
$\nu$ defined by
\begin{equation}
\label{eq:3.8}
I_{\nu} (y) = \sum_{m=0}^{\infty}
\frac{1}{m! \Gamma(m+ \nu +1)} 
\left( 
\frac{y}{2}
\right)^{2m+\nu}.
\end{equation}
By changing the variable $z = ty$
we see that
\begin{equation} \label{eq:3.9}
  (K_1(t)\phi)(x) =
\begin{cases}
& \dfrac{e^{-\frac{t}{2}} t}{2}
           \displaystyle\int_{|y|\le 1}
                I_0 \left({t \sqrt{1-y^2}}/2 \right) \,\phi(x+ty)\,dy
\quad (N=1),
\vspace{5pt}\\ 
  & \dfrac{e^{-\frac{t}{2} } t}{2 \pi}
           \displaystyle\int_{|y|\le 1}
                \frac{\cosh \left( 
                {t \sqrt{1-|y|^2}/2} 
                \right)}
                     {\sqrt{1-|y|^2}}\,
                \phi(x+ty)\,dy \quad (N=2),
                \vspace{5pt}\\
& \dfrac{e^{-\frac{t}{2}} t^{2}}{2}
           \displaystyle\int_{|y|\le1}
                I_1 \left({t \sqrt{1-|y|^2}}/2 \right) \,
               \dfrac{\phi(x+ty)}{\sqrt{1-|y|^{2}}}\,dy,  
               \vspace{5pt}\\
& \qquad\qquad\qquad\qquad\quad +
                \dfrac{e^{-\frac{t}{2}}t}{4 \pi}
                \displaystyle\int_{\mathbb{S}^{2}}
                \phi(x+t\omega)\,d \omega 
                \quad (N=3), 
\end{cases}
\end{equation}
where 
$
S^{N-1}:=
\{\omega\in{\mathbb R}^N, |\omega|=1\}
$
with its surface element $d\omega$.
According to \cite{N},
for $N=3$,
we denote $J_{1}(t)\phi$ and $W_{1}(t)\phi$ by 
\begin{align}
& (K_{1}(t) \phi)(x) =(J_{1}(t)\phi)(x) +(W_{1}(t)\phi)(x), \label{eq:3.10} \\
\label{eq:3.11}
(J_{1}(t)\phi)(x) 
& :=
\frac{e^{-\frac{t}{2}} t^{2}}{2}
           \int_{|y|\le 1}
                I_1 ({t \sqrt{1-|y|^2}}/2 ) \,
               \frac{\phi(x+ty)}{\sqrt{1-|y|^{2}}}\,dy, \\
(W_{1}(t)\phi)(x) & :=
\frac{e^{-\frac{t}{2}}t }{4 \pi}
           \int_{\mathbb{S}^{2}}
\phi(x+t\omega)\,d \omega.
\label{eq:3.12}
\end{align}
Furthermore, it follows from the straightforward calculation (see, e.g., \cite{MN,T01})
that the representation formula for $K_{0}(t) \phi$ is given by
\begin{equation} \label{eq:3.13}
\begin{split}
&
(K_{0}(t)\phi)(x) 
= 
(\partial_{t}
K_{1}(t)\phi)(x) 
+
\frac{1}{2}(K_{1}(t)\phi)(x) \\
&\quad
 =  
\begin{cases}
&  \dfrac{1}{2} (K_{1}(t) \phi)(x)
+ \dfrac{e^{-\frac{t}{2}}t}{4}
           \displaystyle\int_{|y|\le1}
    I_1 \left( t\sqrt{1-y^2}/2 \right) \,\frac{\phi(x+ty)}{\sqrt{1-y^{2}}}
               \,dy, 
               \vspace{5pt}\\ 
& \qquad\qquad\qquad\qquad\qquad
+ e^{-\frac{t}{2}}(\phi(x-t)+\phi(x+t)) \qquad (N=1), 
\vspace{8pt}\\  
& t^{-1} (K_{1}(t) \phi)(x)
+ 
\dfrac{e^{-\frac{t}{2}t } t}{2 \pi}
           \displaystyle\int_{|y|\le 1}
                \sinh \left( 
                t {\sqrt{1-|y|^2}}/2 
                \right)
                  \phi(x+ty) 
                \,dy 
                \vspace{5pt}\\
&
\qquad
 + 
\dfrac{e^{-\frac{t}{2} } t}{2 \pi}
           \displaystyle\int_{|y|\le 1}
                \dfrac{\cosh \left( 
                t{\sqrt{1-|y|^2}}/2 
                \right)}
                     {\sqrt{1-|y|^2}}\,
                 \nabla \phi(x+ty) \cdot y
                \,dy \qquad (N=2),
                \vspace{8pt}\\
&
=  2 t^{-1}
(J_{1}(t)\phi)(x)  + 
(\partial_{t}W_{1}(t)\phi)(x)
\vspace{5pt}\\
&
\qquad
+
\dfrac{e^{-\frac{t}{2}} t^{2}}{2}
           \displaystyle\int_{|y|\le 1}
               \partial_{t}
               \left(
                I_1 ({t \sqrt{1-|y|^2}}/2 ) \,
               \frac{\phi(x+ty)}{\sqrt{1-|y|^{2}}}
               \right)
               \,dy \qquad (N=3). 
\end{cases} 
\end{split} 
\end{equation}
Next we begin with mentioning the $L^{q}$-$L^{r}$ estimates 
for the evolution operators $K_{0}(t)$ and $K_{1}(t)$.  

\begin{Lemma}
{\rm(\cite{HO,MN,Na,N}.)}
\label{Lemma:3.3}
Let $1 \le r \le q \le \infty$ and  $0< \delta < 1/2$.
Assume that $(\phi, \psi) \in L^{r} \cap W^{[N/2], q} \times 
L^{r} \cap L^{q}$.
Then it holds that 
\begin{equation}  
\label{eq:3.14}
\begin{split} 
& \| 
K_{0}(t) \phi 
\|_{L^q} 
\le C 
(1+t)^{-\frac{N}{2}(\frac{1}{r}-\frac{1}{q})}
\|
\phi 
\|_{L^r} 
+ C e^{-\delta t}
\|
\phi
 \|_{W^{[N/2],q}}, 
 \qquad t \ge 0, \\
& \| 
K_{1}(t) \psi 
\|_{L^q} \le C(1+t)^{-\frac{N}{2}(\frac{1}{r}-\frac{1}{q})}
\|
\psi \|_{L^r}
+ C e^{-\delta t}
\|
\psi \|_{L^q}, 
\qquad t \ge 0.
\end{split}
\end{equation}
\end{Lemma}
The $L^{1}$ weighted estimates for the operators $K_{0}(t)$ and $K_{1}(t)$ are well-known. 

\begin{Lemma}
\label{Lemma:3.4}
{\rm(\cite[Proposition~5.2]{T01})}
Let $k\ge0$ and $0<\delta<1/8$.
Then they hold that
\begin{align}  
\label{eq:3.15}  
&
\|K_0(t)\phi \|_{L^1_k}
\le  C \|\phi\|_{W^{[N/2],1}_k} +Ct^{\frac{k}{2}} \|\phi\|_{L^1},
\\
&
\label{eq:3.16}
\|K_1(t)\phi \|_{L^1_k}
\le C \|\phi\|_{L^1_k} 
+
Ct^{\frac{k}{2}} \| \phi \|_{L^1}, 
\end{align}
for all $t\ge0$.
\end{Lemma}
\begin{Remark}
\label{Remark:3.1}
The estimates \eqref{eq:3.15} and \eqref{eq:3.16} are already obtained 
in {\rm\cite{KU}} by another method.
\end{Remark}

From now on, let us introduce the sequences of the heat semi-group, $V_{0}(t)$ and $V_{1}(t)$ as
\begin{align}
&
\label{eq:3.17}
V_0(t):=\frac{1}{2}
\sum_{\ell=0}^{[K/2]} 
\sum_{m=0}^{[K/2]-\ell}
\Phi_{\ell,m}
(-t)^{\ell} (-\Delta)^{2\ell+m}e^{t\Delta},
\\
&
\label{eq:3.18}
V_1(t):=
\sum_{\ell=0}^{[K/2]} 
\sum_{m=0}^{[K/2]-\ell}
\sum_{n=0}^{[K/2]-\ell-m}
\Phi_{\ell,m} \Psi_{n} 
(-t)^{\ell} (-\Delta)^{2\ell+m+n}e^{t\Delta},
\end{align}
for $K \ge 0$,
where $(\Phi_{\ell,m},\Psi_n)$ is given in \eqref{eq:2.3}.
The following lemma states that the operators $K_{0}(t)$ and $K_{1}(t)$ are approximated by  
$V_{0}(t)$ and $V_{1}(t)$, respectively.
This fact plays a crucial role in the proof of Theorem \ref{Theorem:2.2}. 

\begin{Lemma}
 {\rm(\cite[Proposition~4.1]{T01}.)}
\label{Lemma:3.5}
Let 
$K\ge0$,
$1 \le r \le q \le \infty$ and $0< \delta <1/2$.
Assume that $(\phi, \psi) \in L^{r} \cap W^{[N/2], q} \times 
L^{r} \cap L^{q}$.
Then it holds that
\begin{equation} 
\label{eq:3.19} 
\begin{split}
&
\left\| (K_{0}(t)-V_0(t))\phi\right\|_{L^q} 
\le 
Ct^{-\frac{N}{2}
(\frac{1}{r} -\frac{1}{q})-([\frac{K}{2}]+1)} \| \phi \|_{L^r}+ C e^{-\delta t}
\| \phi \|_{W^{[N/2],q}},
\qquad t>0,
\\
&  
\left\| (K_{1}(t)-V_1(t))\phi\right\|_{L^q} 
\le 
Ct^{-\frac{N}{2}(\frac{1}{r} -\frac{1}{q})-([\frac{K}{2}]+1)} \| \phi \|_{L^r}+ C e^{-\delta t}
\| \phi \|_{L^q},  \qquad t>0,
\end{split}
\end{equation}
where $V_0(t)$ and $V_1(t)$ are the operators 
given in \eqref{eq:3.17} and \eqref{eq:3.18}, respectively.
\end{Lemma}

We note that using the notation $K_{0}(t)$ and $K_{1}(t)$, 
the solution of \eqref{eq:1.5} is expressed as
\begin{equation}
\label{eq:3.20}
v(x,t)=(K_0(t)v_0)(x)+(K_1(t)\left(\frac{1}{2}v_0+v_1\right))(x)
\end{equation}
for all $(x,t)\in{\mathbb R}^N\times(0,\infty)$.
As a easy consequence of Lemma \ref{Lemma:3.5} with \eqref{eq:3.20}, 
we have the following:
\begin{Proposition}
{\rm(\cite[Theorem~1.1]{T01})}
\label{Proposition:3.2}
Let $K\ge0$ and $(v_{0}, v_{1}) \in W^{[N/2],1}_{K} \cap W^{[N/2], \infty}\times 
L^{1}_{K} \cap L^{\infty}$.
Then there exists a unique solution $v$ of \eqref{eq:1.5} satisfying
$$
v(t)\in C([0, \infty): L^{1}_{K}) \cap L^{\infty}(0,\infty: L^\infty).
$$
Furthermore, for any $q\in[1,\infty]$,
$$
t^{\frac{N}{2}(1-\frac{1}{q})}\left\| v(t)-V_0(t)v_0-V_1(t)\left(\frac{1}{2}v_0+v_1\right)\right\|_{L^q} =
o(t^{-\frac{K}{2}})
$$
as $t \to \infty$.
\end{Proposition}
Combining Propositions~\ref{Proposition:3.1} and \ref{Proposition:3.2}, 
we can easily see the following.

\begin{Corollary}
\label{Corollary:3.1}
Assume the same assumptions as in Proposition~$\ref{Proposition:3.2}$.
Then, for any $q\in[1,\infty]$,
\begin{equation} 
\label{eq:3.21}
t^{\frac{N}{2}(1-\frac{1}{q})}\| v(t)-U_{\rm lin}(t)\|_{L^q} =o(t^{-\frac{K}{2}})
\end{equation}
as $t \to \infty$, where $U_{\rm lin}(x,t)$ is defined by \eqref{eq:2.7}.
\end{Corollary}
\begin{Remark}
\label{Remark:3.2}
By this corollary we see that the function $U_{\rm lin}$ is an asymptotic expansion to 
the solution of \eqref{eq:1.5}, 
i.e. $U_{\rm lin}$ is a linear approximation of the linear part for the solution $u$ for \eqref{eq:1.1}.
\end{Remark}
%
\subsection{Useful formula}
We prepare the point-wise estimate for the modified Bessel function.
\begin{Lemma}
{\rm(\cite{AS})} 
\label{Lemma:3.6}
For any $\nu\in{\mathbb N}_0$, it holds that
\begin{align} \label{eq:3.22}
\left|I_{\nu} \left(
x
\right)
\right|
\le 
\begin{cases}
& 
C x^{\nu}, \qquad 0<x \le 1, \\
&  C 
x^{-\frac{1}{2}}
  e^{x}, \qquad x \ge 1,
\end{cases}
\end{align}
where $I_\nu$ is given in \eqref{eq:3.8}.
\end{Lemma}
The following lemma is useful to obtain decay estimates.
\begin{Lemma} 
 {\rm(\cite[Lemma~2.5]{T01}.)}
\label{Lemma:3.7} 
Let 
$k \ge 0 $.
Then, for any $c>0$ and $m=1,2,3$,  
there exists a constant $C>0$ such that  
\begin{equation}
\int_{|y|\le1} 
\frac{|y|^{k} e^{-c t |y|^{2}}}
{(1-|y|^{2})^{\frac{m}{4}}
} \,dy
\le C
t^{-\frac{N}{2}-\frac{k}{2} },\qquad t>0.  
\label{eq:3.23}
\end{equation} 
\end{Lemma}
%
\section{Weighted estimates for the linear solutions}
In this section we give weighted $L^1$ and $L^\infty$ estimates for the solution $v$ of \eqref{eq:2.11}.
In other words, 
our new ingredient for the proof of the main result is the following weighted $L^{\infty}$ estimates 
for the evolution operators $K_{0}(t)$ and $K_{1}(t)$, 
which are useful for us to estimate the nonlinear term in the nonlinear approximation. 

Our purpose here is to obtain the following proposition 
on the weighted estimates for the solution $v$ of the linear damped wave equation \eqref{eq:1.5}.

\begin{Proposition}
\label{Proposition:4.1}
Let $K\ge0$.
Assume $(v_{0}, v_{1}) \in W^{[N/2],1}_{K} \cap W^{[N/2]_K, \infty}\times L^{1}_{K} \cap L^{\infty}_K$.
Let $v$ be a unique solution of \eqref{eq:1.5}.
Then, for any $k\in[0,K]$,
\begin{align}
\label{eq:4.1}
&
\|v(t)\|_{L^1_k}\le CE_K(1+t)^{\frac{k}{2}},\qquad t\ge0,\\
&
\|v(t)\|_{L^\infty_k}\le CE_K(1+t)^{\frac{k-N}{2}},\qquad t\ge0,
\label{eq:4.2}
\end{align}
where $E_K$ is given in \eqref{eq:2.14} with $(u_0,u_1)$ replaced by $(v_0,v_1)$, respectively.
\end{Proposition}
To this end, we show the following lemmas, which are the weighted estimates
for the evolution operators $K_0(t)$ and $K_1(t)$.
\begin{Lemma}
\label{Lemma:4.1}
Assume the same conditions as in Lemma~$\ref{Lemma:3.4}$.
Then it holds that
\begin{align}  
\label{eq:4.3}  
&
\|K_0(t)\phi \|_{L^\infty_k}
\le  C \|\phi\|_{W^{[N/2],\infty}_k} +Ct^{\frac{k}{2}} \|\phi\|_{L^\infty},
\\
&
\label{eq:4.4}
\|K_1(t)\phi \|_{L^\infty_k}
\le C \|\phi\|_{L^\infty_k} 
+
Ct^{\frac{k}{2}} \| \phi \|_{L^\infty}, 
\end{align}
for all $t\ge0$.
Furthermore it holds that
\begin{align}  
\label{eq:4.5}  
& \|K_0(t)\phi \|_{L^\infty_k}
\le 
C t^{-\frac{N}{2}} \|\phi\|_{W^{[N/2],1}_k} 
+Ct^{\frac{k-N}{2}} \| \phi \|_{L^1} 
+C e^{-\delta t} \| \phi\|_{L^{\infty}_k},
\\
\label{eq:4.6}
& \|K_1(t)\phi \|_{L^\infty_k} \le 
 C t^{-\frac{N}{2}} \|\phi\|_{L^1_k} 
+Ct^{\frac{k-N}{2}}\|\phi\|_{L^1} 
+ C e^{-\delta t} \|\phi\|_{L^{\infty}_k},
\end{align}
for all $t\ge1$.
\end{Lemma}
By Lemma~\ref{Lemma:4.1}, for the weighted $L^\infty$ estimates, 
we can remove the singularity at $t=0$, and we obtain the following.
\begin{Corollary} 
\label{Corollary:4.1}
Assume the same assumptions as in Lemma~$\ref{Lemma:3.4}$.
Then it holds that
\begin{align} 
&
\| K_0(t)\phi \|_{L^\infty_k}\le
C(1+t)^{-\frac{N}{2}}(\|\phi\|_{W^{[N/2],1}_k}+\|\phi\|_{W^{[N/2],\infty}_k})
\notag\\
&
\qquad\qquad\qquad
 +C(1+t)^{\frac{k-N}{2}}(\|\phi\|_{W^{[N/2],1}}+\|\phi\|_{W^{[N/2],\infty}}),
 \label{eq:4.7}
\\
&
\| K_1(t)\phi \|_{L^\infty_k}\le
C(1+t)^{-\frac{N}{2}}(\|\phi\|_{L^1_k}+\|\phi\|_{L^\infty_k})
\notag\\
&
\qquad\qquad\qquad
 +C(1+t)^{\frac{k-N}{2}}(\|\phi\|_{L^1}+\|\phi\|_{L^\infty}),
 \label{eq:4.8}
\end{align}
for all $t\ge0$.
\end{Corollary}
By Lemma~\ref{Lemma:3.4} and  Corollary~\ref{Corollary:4.1} with the aid of \eqref{eq:3.20} 
we can prove Proposition~\ref{Proposition:4.1} immediately.
So it suffices to prove Lemma~\ref{Lemma:4.1}.
\vspace{5pt}

Here, for the simplicity of the notation, 
we introduce the auxiliary functional $\mathcal{K}^{(m)}[\phi]$ as follows:
\begin{equation} 
\label{eq:4.9}
\mathcal{K}^{(m)}[\phi](x,t):=
 t^{\frac{N}{2}}
\int_{|y|\le 1} 
\frac{
e^{-\frac{t}{2}|y|^{2}}
 }{(1-|y|^{2})^{\frac{m}{4}}} |\phi(x+ty)|\,dy, \qquad (x,t)\in{\mathbb R}^N\times(0,\infty),
\end{equation}
where $m=0,1,2,3$.
Then we have the weighted $L^\infty$ estimates of $\mathcal{K}^{(m)}[\phi](x,t)$.
\begin{Lemma}
\label{Lemma:4.2}
Assume the same conditions as in Lemma~$\ref{Lemma:3.4}$.
Then, for any $m=0,1,2,3$,
\begin{align} 
\label{eq:4.10}
\| \mathcal{K}^{(m)}[\phi](t) \|_{L^\infty_k}
\le C \| \phi \|_{L^\infty_k} 
+
Ct^{\frac{k}{2}} \| \phi \|_{L^\infty}, \qquad t \ge 0, 
\end{align}
and
\begin{equation} \label{eq:4.11}
\| \mathcal{K}^{(m)}[\phi](t) \|_{L^\infty_k}
 \le 
 C t^{-\frac{N}{2}} \| \phi \|_{L^1_k} 
+
Ct^{\frac{k-N}{2}} \| \phi \|_{L^1} 
+ C e^{-\delta t}\| \phi\|_{L^{\infty}_k}, \quad t >0.
\end{equation}
\end{Lemma}
\noindent
{\bf Proof.}
By \eqref{eq:4.9} we see that
\begin{equation}
\label{eq:4.12}
\begin{split}
|x|^k |\mathcal{K}^{(m)}[\phi](x,t)|  
&
\le C
 t^{\frac{N}{2}}
\int_{|y|\le 1} (|x+ty|^k+|ty|^k)
\frac{
e^{-\frac{t}{2}|y|^{2}}
 }{(1-|y|^{2})^{\frac{m}{4}}} |\phi(x+ty)|\,dy
\\
&=: 
CI_1(x,t) + CI_2(x,t)
\end{split}
\end{equation}
for all $(x,t)\in{\mathbb R}^N\times(0,\infty)$.
\vspace{3pt}

We first show the estimate \eqref{eq:4.10}.
Applying the Minkowski inequality with \eqref{eq:3.23} we have
\begin{align} 
\label{eq:4.13}
&
\left\|
I_1(t)
\right\|_{L^\infty}
\le C t^{\frac{N}{2}}
\int_{|y|\le1} 
\frac{
e^{-\frac{t}{2} |y|^{2}}
 }{(1-|y|^{2})^{\frac{m}{4}}} 
\left\||\cdot+ty|^k \phi(\cdot+ty)
\right\|_{L^\infty}\,dy 
 \le C \| \phi \|_{L^\infty_k}, 
 \\
 \label{eq:4.14}
 &
\left\|I_2(t)\right\|_{L^\infty} 
 \le C t^{\frac{N}{2}+k}
\int_{|y|\le1} 
\frac{
e^{-\frac{t}{2} |y|^{2}} |y|^k
 }{(1-|y|^{2})^{\frac{m}{4}}} 
\left\|\phi(\cdot+ty)
\right\|_{L^\infty}\,dy 
 \le C t^{\frac{k}{2}} \| \phi \|_{L^\infty}.
\end{align}
Then combining the above estimates \eqref{eq:4.12}-\eqref{eq:4.14}, 
we obtain the estimate \eqref{eq:4.10}.
\vspace{3pt}

Next we prove the estimate \eqref{eq:4.11}.
Observing the estimate \eqref{eq:4.12}, 
we estimate $I_1$ and $I_2$, respectively. 
For the term $I_1$, we decompose into the two parts:
\begin{equation*}
\begin{split}
I_1(x,t) 
& \le 
C t^{\frac{N}{2}}
\int_{|y|\le1/2} 
|x+ty|^k 
e^{-\frac{t}{2}|y|^{2}}
 |\phi(x+ty)|\,dy \\
&
\qquad
 + C t^{\frac{N}{2}} e^{-\frac{t}{8}}
\int_{1/2\le|y|\le1}(1-|y|^{2})^{-\frac{m}{4}}
|x+ty|^k |\phi(x+ty)|\,dy \\
& =:I_{11}(x,t) + I_{12}(x,t)  
\end{split}
\end{equation*}
for all $(x,t)\in{\mathbb R}^N\times(0,\infty)$.
By changing the integral variable $z=x+ty$
we see that 
\begin{equation} 
\label{eq:4.15}
\begin{split}
I_{11}(x,t) 
=
C t^{-\frac{N}{2}}
\int_{|z-x|\le t/2} 
|z|^k 
e^{-\frac{|z-x|^{2}}{2t}}
 |\phi(z)|\,dz \le
C t^{-\frac{N}{2}} \| \phi \|_{L^1_k}
\end{split}
\end{equation}
for all $(x,t)\in{\mathbb R}^N\times(0,\infty)$.
For $I_{12}(x,t)$, 
since
\begin{align} 
\label{eq:4.16}
\int_{1/2\le|y|\le1}(1-|y|^{2})^{-\frac{m}{4}}\,dy \le C
\end{align}
for any $m =0,1,2,3$,
we have 
\begin{equation} 
\label{eq:4.17}
\begin{split}
I_{12}(x,t) 
& \le 
 C t^{\frac{N}{2}} e^{-\frac{t}{8}} \|\phi\|_{L^{\infty}_k}
\int_{1/2\le|y\le1}(1-|y|^{2})^{-\frac{m}{4}}\,dy \\
& \le 
 C t^{\frac{N}{2}} e^{-\frac{t}{8}} \| \phi \|_{L^{\infty}_k}
\end{split}
\end{equation}
for all $(x,t)\in{\mathbb R}^N\times(0,\infty)$.
$I_{2}(x,t) $ is estimated by the similar way. 
Indeed, we again decompose $I_{2}(x,t) $ into two parts:
\begin{equation*}
\begin{split}
I_2(x,t) 
& \le 
C t^{\frac{N}{2}}
\int_{|y|\le1/2)} 
|ty|^k
e^{-\frac{t}{2} |y|^{2}}
 |\phi(x+ty)|\,dy \\
&
\qquad
 + C t^{\frac{N}{2}} e^{-\frac{t}{8}}
\int_{1/2\le|y|\le1}
(1-|y|^{2})^{-\frac{m}{4}} 
|ty|^k |\phi(x+ty)|\,dy \\
& =:I_{21}(x,t) + I_{22}(x,t)  
\end{split}
\end{equation*}
for all $(x,t)\in{\mathbb R}^N\times(0,\infty)$.
By the well-known estimate
\begin{align*}
|y|^k 
e^{-\frac{t}{2}|y|^{2}}
\le C t^{-\frac{t}{2}}
\end{align*}
and changing the integral variable $z=x+ty$, 
we arrive at the estimate 
\begin{equation} \label{eq:4.18}
\begin{split}
I_{21}(x,t) 
& = 
C t^{\frac{N}{2}+k}
\int_{|y|\le1/2} 
|y|^k 
e^{-\frac{t}{2}|y|^{2}}
 |\phi(x+ty)|\,dy \\
& \le 
C t^{\frac{N+k}{2}}
\int_{|y|\le1/2} 
 |\phi(x+ty)|\,dy \\
& \le 
C t^{\frac{k-N}{2}}
\int_{{\mathbb R}^N} 
 |\phi(z)|\,dz \le C t^{\frac{k-N}{2}} \|\phi\|_{L^1} 
\end{split}
\end{equation}
for all $(x,t)\in{\mathbb R}^N\times(0,\infty)$.
$I_{22}(x,t) $ is estimated easily by the estimate \eqref{eq:4.16} and $|y| \le 1$:
\begin{equation} 
\label{eq:4.19}
I_{22}(x,t) 
\le 
 C t^{\frac{N}{2}} e^{-\frac{t}{8}} \|\phi\|_{L^\infty}
\int_{1/2\le|y|\le1}(1-|y|^{2})^{-\frac{m}{4}}\,dy 
\le C t^{\frac{N}{2}} e^{-\frac{t}{8}} \|\phi\|_{L^\infty}
\end{equation}
for all $(x,t)\in{\mathbb R}^N\times(0,\infty)$.
Then, 
summing up the estimates \eqref{eq:4.12}, \eqref{eq:4.15}, \eqref{eq:4.17}, \eqref{eq:4.18} and \eqref{eq:4.19},
we have 
$$
|x|^k |\mathcal{K}^{(m)}[\phi](x,t)|
\le     
C t^{-\frac{N}{2}} \| \phi \|_{L^1_k}
+C t^{\frac{k-N}{2}} \|\phi\|_{L^1}
+C t^{\frac{N}{2}} e^{-\frac{t}{8}}\|\phi\|_{L^\infty_k},
$$
which implies the estimate \eqref{eq:4.11}. 
Thus the proof of Lemma~\ref{Lemma:4.2} is complete.
$\Box$
\vspace{8pt}

As a next step, we show the point-wise estimates for the evolution operators $K_{1}(t)$ and $K_{0}(t)$ 
by the auxiliary functional $\mathcal{K}^{(m)}[\phi](x,t)$.
\begin{Lemma} 
\label{Lemma:4.3}
Let $\mathcal{K}^m[\phi]$ be functions given in \eqref{eq:4.9}.
Then, for the case $N=1,2$,
\begin{equation} 
\label{eq:4.20}
|(K_{1}(t) \phi)(x) | \le C 
\mathcal{K}^{(N)}[\phi](x,t),\qquad (x,t)\in{\mathbb R}^N\times(0,\infty),
\end{equation}
and, for the case $N=3$,
\begin{equation}
\label{eq:4.21}
| (J_{1}(t) \phi)(x) | \le C \mathcal{K}^{(N)}[\phi](x,t),\qquad (x,t)\in{\mathbb R}^N\times(0,\infty).
\end{equation}
\end{Lemma}
\noindent
{\bf Proof.}
The proof is an easy consequence of Lemma \ref{Lemma:3.6}.
By \eqref{eq:3.22} we see that
\begin{equation} 
\label{eq:4.22}
\left|I_{\nu} \left(
x
\right)
\right|
\le C 
x^{-\frac{1}{2}} e^{x}
\end{equation}
for $\nu \in {\mathbb N}_0$ and $x>0$.
Furthermore, for any $y\in\{y\in{\mathbb R}^N\,:\,|y|<1\}$, 
we have
\begin{equation} 
\label{eq:4.23}
-\frac{1}{2}
+
\frac{ \sqrt{1-|y|^2}}{2}
= 
\frac{-|y|^{2}}
{2(1+ \sqrt{1-|y|^2})} \le -\frac{|y|^{2}}{2}.
\end{equation}
We first prove \eqref{eq:4.20}.
For $N=1$, 
applying the estimate \eqref{eq:4.22} to the expression \eqref{eq:3.9} and \eqref{eq:4.23},
we have 
\begin{equation*}
\begin{split}
|(K_{1}(t)\phi)(x)| 
& \le 
C e^{-\frac{t}{2}} t
\int_{|y|\le1}
\frac{
e^{\frac{t \sqrt{1-y^{2}}}{2}}
}{t^{\frac{1}{2}} (1-y^{2})^{\frac{1}{4}}}
|\phi(x+ty)| \,dy \\
& \le 
C t^{\frac{1}{2}}
\int_{|y|\le1}
\frac{
e^{-\frac{t}{2}|y|^{2}}
}{(1-y^{2})^{\frac{1}{4}}}
|\phi(x+ty)| \,dy
\end{split}
\end{equation*}
for all $(x,t)\in{\mathbb R}^N\times(0,\infty)$.
Therefore we obtain the desired estimate \eqref{eq:4.20} for $N=1$.
Furthermore,
by the definition of $\cosh y$, 
we see that 
\begin{equation} 
\label{eq:4.24}
0 \le \cosh y \le e^{y},\qquad y\in{\mathbb R}.
\end{equation}
Thus, again using the estimate \eqref{eq:4.23}, 
we can easily have the desired estimate \eqref{eq:4.20} for $N=2$.
Next we prove the estimate \eqref{eq:4.21}.
Applying the estimates \eqref{eq:4.22} and \eqref{eq:4.23} to the expression \eqref{eq:3.11}
we have 
\begin{equation*}
\begin{split}
|(J_{1}(t)\phi)(x)| 
& \le 
C e^{-\frac{t}{2}} t^{2}
\int_{|y|\le1}
\frac{
e^{\frac{t \sqrt{1-y^{2}}}{2}}
}{t^{\frac{1}{2}} (1-y^{2})^{\frac{3}{4}}}
|\phi(x+ty)| \,dy \\
& \le 
C t^{\frac{3}{2}}
\int_{|y|\le1}
\frac{
e^{-\frac{1}{2}t|y|^{2}}
}{(1-y^{2})^{\frac{3}{4}}}
|\phi(x+ty)| \,dy
\end{split}
\end{equation*}
for all $(x,t)\in{\mathbb R}^N\times(0,\infty)$,
which is the desired estimate.
Then Lemma~\ref{Lemma:4.3} follows.
$\Box$
\vspace{8pt}
\begin{Lemma}
\label{Lemma:4.4}
Assume the same assumption as in Lemma~$\ref{Lemma:4.3}$.
Then the following holds:
\begin{itemize}
\item[\rm(i)]
For $N=1$,
\begin{equation}
\label{eq:4.25}
|(K_{0}(t)\phi)(x)
\le
C |(K_{1}(t) \phi)(x)|
+ Ce^{-\frac{t}{2}}(|\phi(x-t)|+|\phi(x+t)|)
+ C \mathcal{K}^{(3)}[\phi](x,t) 
\end{equation}
for all $(x,t)\in{\mathbb R}\times(0,\infty)$;
\item[\rm(ii)]
For $N=2$,
\begin{equation}
\label{eq:4.26}
|(K_{0}(t)\phi)(x)
\le
Ct^{-1} |(K_{1}(t) \phi)(x)|
+ C \mathcal{K}^{(0)}[\phi](x,t)
+ C \mathcal{K}^{(2)}[|\nabla \phi|](x,t)
\end{equation}
for all $(x,t)\in{\mathbb R}^2\times(0,\infty)$;
\item[\rm(iii)]
For $N=3$,
\begin{equation}
\label{eq:4.27}
\begin{split}
|(K_{0}(t)\phi)(x)
&
\le
C (t^{-1}
|(J_{1}(t)\phi)(x) | +
 \mathcal{K}^{(1)}(x,t)
+ |(K_{1}(t)|\nabla \phi|)(x)| \\
& \qquad \qquad\qquad\qquad
+
 |(W_{1}(t)\phi|)(x)| 
+
 |(W_{1}(t) |\nabla \phi|)(x) |)
 \end{split}
\end{equation}
for all $(x,t)\in{\mathbb R}^3\times(0,\infty)$.
\end{itemize}
\end{Lemma}
\noindent
{\bf Proof.}
We first prove for the case $N=1$. 
Observing the estimates \eqref{eq:4.10},  
we use the estimates \eqref{eq:3.22} and \eqref{eq:4.23} to see that
\begin{equation} \label{eq:4.28}
\begin{split}
|(K_{0}(t)\phi)(x)| 
& \le
C |(K_{1}(t) \phi)(x)|
+ e^{-\frac{t}{2}}(|\phi(x-t)|+|\phi(x+t)|)
\\
&\qquad
+ Ce^{-\frac{t}{2}} t^{\frac{1}{2}}
           \int_{|y|\le1}
    e^{\frac{t \sqrt{1-y^{2}}}{2}}(1-y^{2})^{-\frac{3}{4}}|\phi(x+ty)|\,dy
\end{split}
\end{equation}
for all $(x,t)\in{\mathbb R}\times(0,\infty)$.
In addition,
by \eqref{eq:4.22} and \eqref{eq:4.23}, it is easy to see that
$$
 Ce^{-\frac{t}{2}} t^{\frac{1}{2}}
           \int_{|y|\le1}
    e^{\frac{t \sqrt{1-y^{2}}}{2}}(1-y^{2})^{-\frac{3}{4}}|\phi(x+ty)|\,dy
    \le C \mathcal{K}^{(3)}(x,t)
$$
for all $(x,t)\in{\mathbb R}\times(0,\infty)$.
This together with \eqref{eq:4.28}  yields the desired estimate \eqref{eq:4.25}.

Next we show the point-wise estimate for $K_{0}(t)g$ with $N=2$.  
We recall the estimate \eqref{eq:4.24} and 
$$
0 \le \sinh y \le e^{y}
$$
for $y \ge 0$.
Then we use \eqref{eq:3.13} to have 
\begin{equation*}
\begin{split}
|(K_{0}(t)\phi)(x)| 
& \le 
t^{-1} |(K_{1}(t)\phi)(x)|
+C e^{-\frac{t}{2}} t
           \displaystyle\int_{|y|\le1}
                e^{ t \frac{\sqrt{1-|y|^2}}{2} }
                  |\phi(x+ty)| 
                \,dy \\
&\qquad
 + Ce^{-\frac{t}{2} } t
           \displaystyle\int_{|y|\le1}
                \dfrac{
e^{ t \frac{\sqrt{1-|y|^2}}{2} }
                     }
                     {\sqrt{1-|y|^2}}\,
                 |\nabla \phi(x+ty)| |y|
                \,dy \\
& \le 
t^{-1} |(K_{1}(t)\phi)(x)|
+ C \mathcal{K}^{(0)}[\phi](x,t)
+ C \mathcal{K}^{(2)}[|\nabla \phi|](x,t),
\end{split}
\end{equation*}
which is the desired estimate \eqref{eq:4.26}.

Finally we show the case $N=3$.
To this end, 
we use the point-wise estimate, which is given in \cite[(5.13)]{T01}, as follows:
\begin{equation*}
\begin{split}
|(K_{0}(t)\phi)(x)| 
\le 
&
C t^{-1}
|(J_{1}(t)\phi)(x) |
+
C
t^{\frac{3}{2}}
           \int_{|y|\le1}
          \frac{
           e^{ \frac{-t |y|^{2}}{2(1+\sqrt{1-|y|^2})} } 
          }{ (\sqrt{1-|y|^2})^{\frac{1}{2}}}
               |\phi(x+ty)|
               \,dy\\
               & 
+ |(K_{1}(t)|\nabla \phi|)(x)|
+
C |(W_{1}(t)\phi)(x)| 
+
C |(W_{1}(t) |\nabla \phi|)(x) |
\end{split}
\end{equation*}
Then by the definition of $\mathcal{K}^{(1)}$, we obtain the estimate \eqref{eq:4.27}.
Thus the proof of Lemma~\ref{Lemma:4.4} is complete.
$\Box$
\vspace{8pt}

Finally we show Lemma~\ref{Lemma:4.1}.
\vspace{5pt}

\noindent
{\bf Proof of Lemma~\ref{Lemma:4.1}.}
We first prove the case $N=1$.
By Lemma~\ref{Lemma:4.2} with \eqref{eq:4.20} we have \eqref{eq:4.4}.
Since
\begin{equation*}
\begin{split}
& |x|^k(|\phi(x-t)|+|\phi(x+t)|) \\
& \le C(|x-t|^k|\phi(x-t)|+t^k |\phi(x-t)|+|x+t|^k|\phi(x+t)| +t^k |\phi(x+t)|),
\end{split}
\end{equation*}
we see that
\begin{equation} 
\label{eq:4.29}
\|(|\phi(\cdot-t)|+|\phi(\cdot+t)|) \|_{L^\infty_k} \le C(|\| \phi \|_{L^\infty_k}+t^k \|\phi \|_{L^\infty}),\qquad t\ge0.
\end{equation}
This together with \eqref{eq:4.4}, \eqref{eq:4.10} and \eqref{eq:4.25} yields 
the estimate \eqref{eq:4.3}.
Furthermore, 
combining the estimates \eqref{eq:4.11}, \eqref{eq:4.20}, \eqref{eq:4.25} and \eqref{eq:4.29}, 
we have the estimates \eqref{eq:4.5} and \eqref{eq:4.6}. 

Next we show the case $N=2$.
Similarly to the case $N=1$,
by Lemma~\ref{Lemma:4.2} with \eqref{eq:4.20} we have \eqref{eq:4.4}.
By \eqref{eq:4.23} and \eqref{eq:4.24}
we observe that 
\begin{equation} 
\label{eq:4.30}
\begin{split}
|(K_{1}(t)\phi)(x)|
& \le Cte^{-\frac{t}{2}} 
           \displaystyle\int_{|y|\le1}
                \dfrac{\cosh \left(t \sqrt{1-|y|^2}/2 
                \right)}
                     {\sqrt{1-|y|^2}}\,
                |\phi(x+ty)| \,dy \\
& \le Ct\displaystyle\int_{|y|\le1}
                \frac{|\phi(x+ty)|}{\sqrt{1-|y|^2}}\, \,dy
\end{split}
\end{equation}
for all $(x,t)\in{\mathbb R}^2\times(0,\infty)$.
Then, by \eqref{eq:3.23} and applying the Minkowski inequality to \eqref{eq:4.30}, we see that 
\begin{equation*} 
\begin{split}
& |x|^k |(K_{1}(t)\phi)(x)| \\
& \le Ct \int_{|y|\le1}
\frac{
|x+ty|^k |\phi(x+ty) | +|ty|^k |\phi(x+ty)| }{\sqrt{1-|y|^2}} \,dy \\
& \le Ct \|\phi\|_{L^\infty_k} + C t^{\frac{k}{2}} \|\phi \|_{L^\infty}
\end{split}
\end{equation*}
for all $(x,t)\in{\mathbb R}^2\times(0,\infty)$.
This together with \eqref{eq:4.10} and \eqref{eq:4.26} implies \eqref{eq:4.3}.  
Furthermore, 
combining the estimates \eqref{eq:4.11}, \eqref{eq:4.20} and \eqref{eq:4.26}, 
we have the estimates \eqref{eq:4.5} and \eqref{eq:4.6}. 

Finally we prove the case $N=3$.
For the term $W_{1}(t) g$, 
by \eqref{eq:3.12}
we have
\begin{equation} 
\label{eq:4.31}
\begin{split}
 |x|^k |(W_{1}(t)\phi)(x)| 
& \le 
C 
e^{-\frac{t}{2} } t 
\int_{\mathbb{S}^{2}}
(|x+t\omega|^k +t^k)
|\phi(x+t\omega)|\,d \omega \\
& \le 
C 
e^{-\frac{t}{2} } t
\left( 
\|\phi\|_{L^\infty_k}
+
t^k \|\phi\|_{L^\infty}
\right)
\end{split}
\end{equation}
 for all $(x,t)\in{\mathbb R}^3\times(0,\infty)$.
Then, noting \eqref{eq:3.10}, \eqref{eq:4.10}, \eqref{eq:4.21} and \eqref{eq:4.31}, 
we arrive at the estimate 
\begin{equation*}
\begin{split}
\| K_{1}(t) \phi \|_{L^\infty_k} 
& \le C 
\|\mathcal{K}^{(3)}[\phi](t) \|_{L^\infty_k}
+C 
e^{-\frac{t}{2} } t
\left( 
\|\phi\|_{L^\infty_k}
+
t^k \|\phi\|_{L^\infty}
\right) \\
& \le  C \|\phi\|_{L^\infty_k} 
+
Ct^{\frac{k}{2}} \|\phi\|_{L^\infty}, \qquad t \ge 0, 
\end{split}
\end{equation*} 
which is the estimate \eqref{eq:4.4}.
To show the estimate \eqref{eq:4.3}, 
we recall that the estimate \eqref{eq:3.22} implies 
\begin{align*} 
\left|I_{\nu} \left(
x
\right)
\right|
\le Ce^{x},\qquad x>0,
\end{align*}
for $\nu \in \mathbb{N}_0$.
Then, by \eqref{eq:3.23} we obtain
\begin{equation} 
\label{eq:4.32}
\begin{split}
& |x|^k |(J_{1}(t)\phi)(x)| \\
& \le Ce^{-\frac{t}{2}}t
           \int_{|y|\le1} e^{\frac{t\sqrt{1-|y|^{2}}}{2}}\,
               \frac{|x|^k |\phi(x+ty)|}{\sqrt{1-|y|^{2}}}\,dy \\
& \le Ct \int_{|y|\le1}
\frac{
(|x+ty|^k+|ty|^k) |\phi(x+ty)|}{\sqrt{1-|y|^2}}\, \,dy
\le Ct( \|\phi\|_{L^\infty_k} + t^{\frac{k}{2}} \|\phi\|_{L^\infty})
\end{split}
\end{equation}
for all $(x,t)\in{\mathbb R}^3\times(0,\infty)$.
Summing up the estimates \eqref{eq:4.10}, \eqref{eq:4.20}, \eqref{eq:4.27}, 
\eqref{eq:4.31} and \eqref{eq:4.32},
we obtain the estimate \eqref{eq:4.3}. 
By \eqref{eq:3.10}, \eqref{eq:4.11}, \eqref{eq:4.21} and \eqref{eq:4.31}
we reached at the estimate \eqref{eq:4.6}.
Furthermore,
by \eqref{eq:4.6}, \eqref{eq:4.11}, \eqref{eq:4.20} and \eqref{eq:4.31}
we can easily obtain the estimate \eqref{eq:4.5}.
Thus Lemma~\ref{Lemma:4.1} follows.
$\Box$
\vspace{8pt}
%
\section{Proof of Theorem \ref{Theorem:2.1}}
This section is devoted to the proof of Theorem \ref{Theorem:2.1}.
We note that the our strategy for the proof of the estimates \eqref{eq:2.12} and \eqref{eq:2.13}
are the combination of the  argument similar to the proof of \cite[Lemma~4.1]{KU} and 
the weighted $L^{\infty}$ estimates for the linearized solution which is developed in Section 4.  
\vspace{5pt}

First we give the following weighted $L^1$ estimate for the solution to \eqref{eq:1.1}.

\begin{lemma}
\label{Lemma:5.1}
Assume the same conditions as in Theorem $\ref{Theorem:2.1}$.
Then there exists a positive constant $\varepsilon_0$ such that,
if $E_0\le\varepsilon_0$,
then the estimate \eqref{eq:2.12} holds.
\end{lemma}
\noindent
{\bf Proof.}
For any $t \ge 0$, we put
$$
E(t) := \sup_{0\le s \le t}(1+s)^{-\frac{k}{2}}\|u(s)\|_{L^1_k}.
$$
Then, by \eqref{eq:2.4}, \eqref{eq:4.1} and \eqref{eq:3.16} we obtain
\begin{equation}
\label{eq:5.1}
\|u(t)\|_{L^1_k}\le CE_K(1+t)^{\frac{k}{2}}
+C \int_0^t (\|F(u(s))\|_{L^1_k} + (t-s)^{\frac{k}{2}}\|F(u(s))\|_{L^1})\,ds
\end{equation}
for all $t\ge0$.
Here, by \eqref{eq:1.2} and \eqref{eq:1.3}
we have
\begin{align}
&
\|F(u(t))\|_{L^1} \le C \|u(t)\|_{L^p}^p
\le CE_0^p(1+t)^{-A}, 
\label{eq:5.2} \\[2mm]
&
\|F(u(t))\|_{L^1_k} \le C\|u(t)\|_{L^\infty}^{p-1}\|u(t)\|_{L^1_k}
\le C E_0^{p-1} 
(1+t)^{-A} \|u(t)\|_{L^1_k},
\label{eq:5.3} 
\end{align}
for all $t\ge0$, where $A$ is defined by \eqref{eq:2.15}.
Since $A>1$,
substituting the above estimates into \eqref{eq:5.1}, this yields
\begin{equation*}
\begin{split}
\|u(t)\|_{L^1_k}
&
\le C E_K(1+t)^{\frac{k}{2}}
+ C E_0^{p-1} E(t) \int_0^t (1+s)^{\frac{k}{2}-A}\,ds
+ C E_0^p \int_0^t (t-s)^{\frac{k}{2}}(1+s)^{-A}\,ds \\
&\le C(E_K + E_0^p+E_0^{p-1} E(t))(1+t)^{\frac{k}{2}}, \quad t\ge0. 
\end{split}
\end{equation*}
Eventually, we obtain 
$$
E(t) \le C_0(E_K + E_0^p+E_0^{p-1}E(t)),\qquad t\ge0,
$$
where $C_0$ is a positive constant.
Choosing $\varepsilon_0$ such that
$$
0<\varepsilon_0^{p-1}\le\min\{1/(2C_0),E_K^{(p-1)/p)}\},
$$ 
and letting $E_0$ sufficiently small such that $E_0 \le \varepsilon_0$, 
we obtain $E(t) \le 4C_0E_K$ and the desired estimate \eqref{eq:2.12}.
Thus Lemma \ref{Lemma:5.1} follows.
$\Box$
\vspace{8pt}

Secondly, we give the following weighted $L^\infty$ estimate for the solution $u$ of \eqref{eq:1.1}.
\begin{lemma}
\label{Lemma:5.2}
Assume the same conditions as in Theorem $\ref{Theorem:2.1}$.
Then there exists a positive constant $\varepsilon_1$ such that,
if $E_0\le\varepsilon_1$,
then the estimate \eqref{eq:2.13} holds.
\end{lemma}
\noindent
{\bf Proof.}
For any $t \ge 0$, we put
$$
\tilde{E}(t) := \sup_{0\le s \le t}(1+s)^{\frac{N-k}{2}}\|u(s)\|_{L^\infty_k}.
$$
Then, by \eqref{eq:2.4} and \eqref{eq:4.2} we obtain
\begin{equation}
\label{eq:5.4}
\|u(t)\|_{L^\infty_k}
\le CE_K(1+t)^{\frac{k-N}{2}}
+\left(\int_0^{t/2}+\int_{t/2}^t\right) \|K_1(t-s)F(u(s))\|_{L^\infty_k}\,ds
\end{equation}
for all $t\ge0$.
Here, by \eqref{eq:1.2} and \eqref{eq:1.3}
we have
\begin{align}
&
\|F(u(t))\|_{L^\infty} \le C \|u(t)\|_{L^\infty}^p
\le CE_0^p(1+t)^{-A-1-\frac{N}{2}}, 
\label{eq:5.5} \\[2mm]
&
\|F(u(t))\|_{L^\infty_k} \le C\|u(t)\|_{L^\infty}^{p-1}\|u(t)\|_{L^\infty_k}
\le C E_0^{p-1} 
(1+t)^{-A-1} \|u(t)\|_{L^\infty_k},
\label{eq:5.6} 
\end{align}
for all $t\ge0$, where $A$ is defined by \eqref{eq:2.15}.
Then, since $A>0$, by \eqref{eq:4.4}, \eqref{eq:5.5} and \eqref{eq:5.6} we obtain
\begin{equation}
\label{eq:5.7}
\begin{split}
&
\int_{t/2}^t\|K_1(t-s)F(u(s))\|_{L^\infty_k}\,ds
\\
&
\le C\int_{t/2}^t(\|F(u(s))\|_{L^\infty_k}+(t-s)^{\frac{k}{2}}\|F(u(s))\|_{L^\infty})\,ds
\\
&
\le CE_0^{p-1}\tilde{E}(t)\int_{t/2}^t(1+s)^{-A-1-\frac{N-k}{2}}\,ds
+CE_0^p\int_{t/2}^t(t-s)^{\frac{k}{2}}(1+s)^{-A-1-\frac{N}{2}}\,ds
\\
&
\le CE_0^{p-1}(\tilde{E}(t)+E_0)(1+t)^{-\frac{N-k}{2}}
\end{split}
\end{equation}
for all $t\ge0$.
Furthermore,
since $A>0$ again,
by \eqref{eq:2.12}, \eqref{eq:4.8}, \eqref{eq:5.2}, \eqref{eq:5.3}, \eqref{eq:5.5} and \eqref{eq:5.6}
we have 
\begin{equation*}
\begin{split}
&
\int_0^{t/2}\|K_1(t-s)F(u(s))\|_{L^\infty_k}\,ds
\\
&
\le C\int_0^{t/2}(1+t-s)^{-\frac{N}{2}}(\|F(u(s))\|_{L^1_k}+\|F(u(s))\|_{L^\infty_k})\,ds
\\
&
\qquad\qquad
+C\int_0^{t/2}(1+t-s)^{-\frac{N-k}{2}}(\|F(u(s))\|_{L^1}+\|F(u(s))\|_{L^\infty})\,ds
\\
&
\le CE_0^{p-1}(1+t)^{-\frac{N}{2}}\int_0^{t/2}(1+s)^{-A-1}(\|u(s)\|_{L^1_k}+\|u(s)\|_{L^\infty_k}
)\,ds
\\
&
\qquad\qquad
+CE_0^p(1+t)^{-\frac{N-k}{2}}\int_0^{t/2}(1+s)^{-A-1}(1+(1+s)^{-\frac{N}{2}})\,ds
\\
&
\le CE_0^{p-1}E_K(1+t)^{-\frac{N}{2}}\int_0^{t/2}(1+s)^{-A-1+\frac{k}{2}}\,ds
\\
&\qquad\qquad
+CE_0^{p-1}\tilde{E}(t)(1+t)^{-\frac{N}{2}}\int_0^{t/2}(1+s)^{-A-1-\frac{N-k}{2}}\,ds
+CE_0^p(1+t)^{-\frac{N-k}{2}}
\\
&
\le CE_0^{p-1}(E_K+\tilde{E}(t)+E_0)(1+t)^{-\frac{N-k}{2}}
\end{split}
\end{equation*}
for all $t\ge0$.
This together with \eqref{eq:5.4} and \eqref{eq:5.7} implies that
\begin{equation*}
\tilde{E}(t)
\le C_1(E_K + E_0^p+E_0^{p-1}E_K+ E_0^{p-1}\tilde{E}(t)), \quad t\ge0. 
\end{equation*}
where $C_1$ is a positive constant.
Choosing $\varepsilon_1$ such that
$$
0<\varepsilon_1^{p-1}\le\min\{1,1/(2C_1),E_K^{p-1}\}
$$ 
and letting $E_0$ sufficiently small such that $E_0 \le \varepsilon_1$, 
we obtain $\tilde{E}(t) \le 6C_1E_K$ and the desired estimate \eqref{eq:2.13},
and the proof of Lemma \ref{Lemma:5.1} is complete.
$\Box$
\vspace{8pt}

\noindent
{\bf Proof of Theorem~\ref{Theorem:2.1}.}
We note that under the assumption in Theorem \ref{Theorem:2.1}, 
we easily have the existence of global solution to \eqref{eq:1.1} 
satisfying \eqref{eq:2.11} (see e.g. \cite{HO},  \cite{MN}, \cite{N}). 
Therefore, to show Theorem \ref{Theorem:2.1}, 
it suffices to prove the estimates \eqref{eq:2.12} and \eqref{eq:2.13}.
Moreover, by Lemmas \ref{Lemma:5.1} and \ref{Lemma:5.2}, we obtain the estimates \eqref{eq:2.12} and \eqref{eq:2.13} with putting 
$\varepsilon:= \min \{ \varepsilon_{0}, \varepsilon_{1} \}$.
The proof of Theorem \ref{Theorem:2.1} is completed.

\section{Proof of Theorem \ref{Theorem:2.2}}
In this section, applying the argument as in \cite{IKK}, we prove Theorem~\ref{Theorem:2.2}.
The proof consists of two parts, the derivation of the estimates 
\eqref{eq:2.16} and \eqref{eq:2.17}. 
\vspace{5pt}

\noindent
{\bf Proof of (\ref{eq:2.16}).}
The proof is by induction.
Applying estimates \eqref{eq:3.2} and \eqref{eq:3.3}
to the definition of $U_{\rm lin}$, that is \eqref{eq:2.7},
we can easily obtain 
\begin{equation}
\label{eq:6.1}
\sup_{t>0}\,(1+t)^{\frac{N}{2}(1-\frac{1}{q})}\|U_{\rm lin}(t)\|_{L^q}
+\sup_{t>0}\,(1+t)^{-\frac{k}{2}+\frac{N}{2}(1-\frac{1}{\gamma})}\|U_{\rm lin}(t)\|_{L^\gamma_k}
<\infty
\end{equation}
for any $q\in[1,\infty]$, $k\in[0,K]$ and $\gamma\in\{1,\infty\}$.
On the other hand, 
by \eqref{eq:1.2}, \eqref{eq:2.11}, \eqref{eq:2.12}, \eqref{eq:2.13} and \eqref{eq:3.5}
we have
\begin{equation}
\label{eq:6.2}
E_{K,q}[F(u)](t)
\le C\|u(t)\|_{L^\infty}^{p-1}E_{K,q}[u](t)
\le
C(1+t)^{\frac{K}{2}-A-1},
\qquad t>0,
\end{equation}
where $A$ is given in \eqref{eq:2.15}.
By \eqref{eq:3.6} and \eqref{eq:6.2} we obtain
\begin{equation}
\label{eq:6.3}
\begin{split}
&
\left\|\int_0^{t/2}
K_1(t-s)M_\alpha(F(u(s)),s)g_\alpha(s)\,ds\right\|_{L^\gamma_k}
\\
&
\le
\int_0^{t/2}(|M_\alpha F(u(s)),s)|
\|K_1(t-s)g_\alpha(s)\|_{L^\gamma_k}\,ds
\\
&
\le
C\int_0^{t/2}(1+s)^{-\frac{K-|\alpha|}{2}}E_{K,q}[F(u)](s)
\|K_1(t-s)g_\alpha(s)\|_{L^\gamma_k}\,ds
\\
&
\le
C\int_0^{t/2}(1+s)^{-A-1+\frac{|\alpha|}{2}}\|K_1(t-s)g_\alpha(s)\|_{L^\gamma_k}\,ds,
\qquad t>0.
\end{split}
\end{equation}
Here, by \eqref{eq:3.2}, \eqref{eq:3.3}, \eqref{eq:3.16} and \eqref{eq:4.8} we have
\begin{equation}
\label{eq:6.4}
\begin{split}
\|K_1(t-s)g_\alpha(s)\|_{L^1_k}
&
\le C\|g_\alpha(s)\|_{L^1_k}+C(t-s)^{\frac{k}{2}}\|g_\alpha(s)\|_{L^1}
\\
&
\le C(1+t)^{\frac{k-|\alpha|}{2}}
+C(t-s)^{\frac{k}{2}}(1+t)^{-\frac{|\alpha|}{2}},\qquad t>0,
\end{split}
\end{equation}
and
\begin{equation*}
\begin{split}
&
\|K_1(t-s)g_\alpha(s)\|_{L^\infty_k}
\\
&
\le C(1+t-s)^{-\frac{N}{2}}(\|g_\alpha(s)\|_{L^1_k}+\|g_\alpha(s)\|_{L^\infty_k})
\\
&
\qquad\qquad\quad
+C(1+t-s)^{\frac{k-N}{2}}(\|g_\alpha(s)\|_{L^1}+\|g_\alpha(s)\|_{L^\infty})
\\
&
\le C(1+t-s)^{-\frac{N}{2}}(1+s)^{\frac{k-|\alpha|}{2}}
+C(1+t-s)^{\frac{k-N}{2}}(1+s)^{-\frac{|\alpha|}{2}},\qquad t>0.
\end{split}
\end{equation*}
These together with \eqref{eq:6.3} and $A>0$ yield
\begin{equation}
\label{eq:6.5}
\begin{split}
&
\left\|\int_0^{t/2}
K_1(t-s)M_\alpha(F(u(s)),s)g_\alpha(s)\,ds\right\|_{L^\gamma_k}
\\
&
\le
C\int_0^{t/2}
\bigg((1+t-s)^{-\frac{N}{2}(1-\frac{1}{\gamma})}(1+s)^{\frac{k}{2}-A-1}
\\
&
\qquad\qquad\qquad\qquad\qquad
+(1+t-s)^{\frac{k}{2}-\frac{N}{2}(1-\frac{1}{\gamma})}(1+s)^{-A-1}\bigg)\,ds
\\
&
\le
C(1+t)^{\frac{k}{2}-\frac{N}{2}(1-\frac{1}{\gamma})}\int_0^\infty(1+s)^{-A-1}\,ds
\le
C(1+t)^{\frac{k}{2}-\frac{N}{2}(1-\frac{1}{\gamma})},\qquad
t>0.
\end{split}
\end{equation}
Furthermore, 
since it follows from \eqref{eq:3.2}, \eqref{eq:3.3} and \eqref{eq:4.4} that
\begin{equation}
\label{eq:6.6}
\begin{split}
\|K_1(t-s)g_\alpha(s)\|_{L^\infty_k}
&
\le C\|g_\alpha(s)\|_{L^\infty_k}+C(t-s)^{\frac{k}{2}}\|g_\alpha(s)\|_{L^\infty}
\\
&
\le C(1+s)^{\frac{k-N-|\alpha|}{2}}
+C(t-s)^{\frac{k}{2}}(1+s)^{-\frac{N+|\alpha|}{2}},\qquad t>0,
\end{split}
\end{equation}
by \eqref{eq:6.3}, \eqref{eq:6.4}, \eqref{eq:6.6} and $A>0$ we have
\begin{equation*}
\begin{split}
&
\left\|\int_{t/2}^t
K_1(t-s)M_\alpha(F(u(s)),s)g_\alpha(s)\,ds\right\|_{L^\gamma_k}
\\
&
\le
C\int_{t/2}^t
(1+s)^{-A+\frac{|\alpha|}{2}}\|K_1(t-s)g_\alpha(s)\|_{L^\gamma_k}\,ds
\\
&
\le C\int_{t/2}^t
\bigg((1+s)^{-A-1+\frac{k}{2}-\frac{N}{2}(1-\frac{1}{\gamma})}
+(t-s)^{\frac{k}{2}}(1+s)^{-A-1-\frac{N}{2}(1-\frac{1}{\gamma})}\bigg)\,ds
\\
&
\le C(1+t)^{\frac{k}{2}-\frac{N}{2}(1-\frac{1}{\gamma})}
\int_0^\infty((1+s)^{-A-1}\,ds
\le C(1+t)^{\frac{k}{2}-\frac{N}{2}(1-\frac{1}{\gamma})}
\qquad t>0.
\end{split}
\end{equation*}
This together with \eqref{eq:2.9}, \eqref{eq:6.1} and \eqref{eq:6.5} yields
\eqref{eq:2.16} for $j=0$.

Assume that \eqref{eq:2.16} holds for some $j=m\in{\mathbb N}_0$.
Similarly to \eqref{eq:6.2} with \eqref{eq:2.16} for the case $j=m$, we have
\begin{equation}
\label{eq:6.7}
E_{K,q}[F_m](t)
\le C\|U_m(t)\|_{L^\infty}^{p-1}E_{K,q}[U_m](t)
\le
C(1+t)^{\frac{K}{2}-A-1},
\qquad t>0.
\end{equation}
By \eqref{eq:2.9} and \eqref{eq:2.16} with $j=0$ we see that
\begin{equation}
\label{eq:6.8}
\begin{split}
\|U_{m+1}(t)\|_{L^\gamma_k}
&
\le
C(1+t)^{\frac{k}{2}-\frac{N}{2}(1-\frac{1}{r})}
+\left\|\left(\int_0^{t/2}+\int_{t/2}^t\right)K_1(t-s)P_K(s)F_{m}(s)\,ds\right\|_{L^\gamma_k}
\\
&
=:C(1+t)^{\frac{k}{2}-\frac{N}{2}(1-\frac{1}{\gamma})}+J_1(t)+J_2(t),
\qquad t>0.
\end{split}
\end{equation}
For the case $\gamma=1$,
by \eqref{eq:3.7} and \eqref{eq:3.16} we have
\begin{equation}
\label{eq:6.9}
\begin{split}
J_1(t)
&
\le\int_0^{t/2}\|K_1(t-s)P_K(s)F_{m}(s)\|_{L^1_k}\,ds
\\
&
\le C\int_0^{t/2}(\|P_K(s)F_m(s)\|_{L^1_k}+(t-s)^{\frac{k}{2}}\|P_K(s)F_m(s)\|_{L^1})\,ds
\\
&
\le C\int_0^{t/2}((1+s)^{-\frac{K-k}{2}}
+t^{\frac{k}{2}}(1+s)^{-\frac{K}{2}})E_{K,q}[F_m](s)\,ds
\\
&
\le C(1+t)^{\frac{k}{2}}
\int_0^\infty(1+s)^{-\frac{K}{2}}E_{K,q}[F_m](s)\,ds,\qquad t>0.
\end{split}
\end{equation}
Furthermore, for the case $\gamma=\infty$,
by \eqref{eq:3.7} and \eqref{eq:4.8} we obtain
\begin{equation}
\label{eq:6.10}
\begin{split}
J_1(t)
&
\le\int_0^{t/2}\|K_1(t-s)P_K(s)F_{m}(s)\|_{L^\infty_k}\,ds
\\
&
\le C\int_0^{t/2}((1+t-s)^{-\frac{N}{2}}(\|P_K(s)F_m(s)\|_{L^1_k}+\|P_K(s)F_m(s)\|_{L^\infty_k})\,ds
\\
&
\qquad\qquad
+C\int_0^{t/2}(1+t-s)^{\frac{k-N}{2}}(\|P_K(s)F_m(s)\|_{L^1}+\|P_K(s)F_m(s)\|_{L^\infty})\,ds
\\
&
\le C(1+t)^{-\frac{N}{2}}\int_0^{t/2}(1+s)^{-\frac{K-k}{2}}E_{K,q}[F_m](s)\,ds
\\
&
\qquad\qquad\qquad
+C(1+t)^{\frac{k-N}{2}}\int_0^{t/2}(1+s)^{-\frac{K}{2}}E_{K,q}[F_m](s)\,ds
\\
&
\le C(1+t)^{\frac{k-N}{2}}
\int_0^\infty(1+s)^{-\frac{K}{2}}E_{K,q}[F_m](s)\,ds,\qquad t>0.
\end{split}
\end{equation}
Moreover, by \eqref{eq:3.7}, \eqref{eq:3.16} and \eqref{eq:4.4} we have
\begin{equation*}
\begin{split}
J_2(t)
&
\le\int_{t/2}^t\|K_1(t-s)P_K(s)F_{m}(s)\|_{L^\gamma_k}\,ds
\\
&
\le C\int_{t/2}^t\left(\|P_K(s)F_m(s)\|_{L^\gamma_k}
+(t-s)^{\frac{k}{2}}\|P_K(s)F_m(s)\|_{L^\gamma}\right)\,ds
\\
&
\le C\int_{t/2}^t\left((1+s)^{-\frac{K-k}{2}-\frac{N}{2}(1-\frac{1}{\gamma})}
+t^{\frac{k}{2}}(1+s)^{-\frac{K}{2}-\frac{N}{2}(1-\frac{1}{\gamma})}\right)E_{K,q}[F_m](s)\,ds
\\
&
\le C(1+t)^{\frac{k}{2}-\frac{N}{2}(1-\frac{1}{\gamma})}
\int_0^\infty(1+s)^{-\frac{K}{2}}E_{K,q}[F_m](s)\,ds,\qquad t>0.
\end{split}
\end{equation*}
This together with \eqref{eq:6.7}, \eqref{eq:6.9}, \eqref{eq:6.10} and $A>0$ implies that
\begin{equation}
\label{eq:6.11}
J_1(t)+J_2(t)\le C(1+t)^{\frac{k}{2}-\frac{N}{2}(1-\frac{1}{\gamma})}\int_0^\infty(1+s)^{-A-1}\,ds
\le C(1+t)^{\frac{k}{2}-\frac{N}{2}(1-\frac{1}{\gamma})},\qquad t>0.
\end{equation}
Therefore, by \eqref{eq:6.8} and \eqref{eq:6.11}
we obtain \eqref{eq:2.16} with $j=m+1$.
Hence, by induction we see that \eqref{eq:2.16} holds for all $j\in{\mathbb N}_0$.
\vspace{5pt}

\noindent
{\bf Proof of (\ref{eq:2.17}).}
The proof is also by induction.
We put
\begin{equation}
\label{eq:6.12}
u_{\rm lin}(x,t):=(K_0(t)u_0)(x)+\bigg(K_1\bigg(\frac{1}{2}u_0+u_1\bigg)\bigg)(x),
\end{equation}
for all $(x,t)\in{\mathbb R}^N\times(0,\infty)$.
Let $j\in{\mathbb N}_0$ and $F_{-1}\equiv0$.
Since it follows from \eqref{eq:2.4}, \eqref{eq:2.9} and \eqref{eq:6.12} that
\begin{equation*}
\begin{split}
&
u(x,t)-U_j(x,t)
\\
&
=u_{\rm lin}(x,t)+\int_0^t(K_1(t-s)F(u(s)))(x)\,ds
\\
&
\qquad\qquad
-U_0(x,t)-\int_0^t(K_1(t-s)P_K(s)F_{j-1}(s))(x)\,ds
\\
&
=u_{\rm lin}(x,t)-U_{\rm lin}(x,t)+\int_0^t(K_1(t-s)F(u(s)))(x)\,ds
\\
&
\qquad\qquad
-\sum_{|\alpha|\le K}\int_0^t(K_1(t-s)M_\alpha(F(u(s)),s)g_\alpha(s))(s)\,ds
\\
&
\qquad\qquad\qquad\qquad
-\int_0^t(K_1(t-s)P_K(s)F_{j-1}(s))(x)\,ds
\\
&
=u_{\rm lin}(x,t)-U_{\rm lin}(x,t)+\int_0^t\bigg(K_1(t-s)P_K(s)\bigg(F(u(s))-F_{j-1}(s)\bigg)\bigg)(x)\,ds
\end{split}
\end{equation*}
for all $(x,t)\in{\mathbb R}^N\times(0,\infty)$,
for any $q\in[1,\infty]$, we have
\begin{equation}
\label{eq:6.13}
\begin{split}
\|u(t)-U_j(t)\|_{L^q}
&
\le \|u_{\rm lin}(t)-U_{\rm lin}(t)\|_{L^q}
+\left\|\int_0^tK_1(t-s)P_K\tilde{F}_{j-1}(s)\,ds\right\|_{L^q}
\\
&
\le \|u_{\rm lin}(t)-U_{\rm lin}(t)\|_{L^q}
+\left\|\int_{t/2}^tK_1(t-s)P_K\tilde{F}_{j-1}(s)\,ds\right\|_{L^q}
\\
&\qquad\qquad
+\left\|\int_0^{t/2}\bigg(K_1(t-s)-V_1(t-s)\bigg)P_K\tilde{F}_{j-1}(s)\,ds\right\|_{L^q}
\\
&\qquad\qquad\qquad\qquad\qquad\qquad
+\left\|\int_0^{t/2}V_1(t-s)P_K\tilde{F}_{j-1}(s)\,ds\right\|_{L^q}
\\
&
=: \tilde{J}_1(t)+\tilde{J}_2(t)+\tilde{J}_3(t)+\tilde{J_4}(t),
\qquad t>0.
\end{split}
\end{equation}
where 
\begin{equation}
\label{eq:6.14}
\tilde{F}_{j-1}(x,t):=F(u(x,t))-F_{j-1}(x,t).
\end{equation}
By \eqref{eq:3.21} we obtain
\begin{equation}
\label{eq:6.15}
t^{\frac{N}{2}(1-\frac{1}{q})}\tilde{J}_1(t)=o(t^{-\frac{K}{2}})\quad\mbox{as}\quad t\to\infty.
\end{equation}
So it suffices to consider the terms $\tilde{J}_n(t)$ with $n=2,3,4$.

We first consider the case $j=0$.
By \eqref{eq:6.2} and \eqref{eq:6.14} we have
\begin{equation}
\label{eq:6.16}
E_{K,q}[\tilde{F}_{-1}](t)=E_{K,q}[F(u)](t)\le C(1+t)^{\frac{K}{2}-A-1},\qquad t>0.
\end{equation}
Since $A>0$,
it follows from \eqref{eq:3.7}, \eqref{eq:3.14} and \eqref{eq:6.16} that
\begin{equation}
\label{eq:6.17}
\begin{split}
t^{\frac{N}{2}(1-\frac{1}{q})}\tilde{J_2}(t)
&
\le Ct^{\frac{N}{2}(1-\frac{1}{q})}\int_{t/2}^t(1+e^{-\delta(t-s)})\|P_K(s)\tilde{F}_{-1}(s)\|_{L^q}\,ds
\\
&
\le Ct^{\frac{N}{2}(1-\frac{1}{q})}
\int_{t/2}^t(1+s)^{-\frac{K}{2}-\frac{N}{2}(1-\frac{1}{q})}E_{K,q}[\tilde{F}_{-1}](s)\,ds
\\
&
\le C\int_{t/2}^\infty(1+s)^{-A-1}\,ds
\le C(1+t)^{-A},
\qquad t>0.
\end{split}
\end{equation}
Furthermore, 
since $A>0$,
by \eqref{eq:3.7}, \eqref{eq:3.19} and \eqref{eq:6.16} we obtain
\begin{equation}
\label{eq:6.18}
\begin{split}
t^{\frac{N}{2}(1-\frac{1}{q})}\tilde{J_3}(t)
&
\le Ct^{\frac{N}{2}(1-\frac{1}{q})}
\int_0^{t/2}\bigg((t-s)^{-\frac{N}{2}(1-\frac{1}{q})-([\frac{K}{2}]+1)}\|P_K\tilde{F}_{-1}(s)\|_{L^1}
\\
&
\qquad\qquad\qquad\qquad\qquad\qquad\qquad\qquad
+e^{-\delta(t-s)}\|P_K\tilde{F}_{-1}(s)\|_{L^q}\bigg)\,ds
\\
&
\le Ct^{-([\frac{K}{2}]+1)}
\int_0^{t/2}(1+s)^{-\frac{K}{2}}E_{K,q}[\tilde{F}_{-1}](s)\,ds
\\
&
\le Ct^{-([\frac{K}{2}]+1)}\int_0^\infty(1+s)^{-A-1}\,ds
\le Ct^{-([\frac{K}{2}]+1)},
\qquad t>0.
\end{split}
\end{equation}
On the other hand,
it follows from \eqref{eq:3.18} that
\begin{equation*}
\begin{split}
\|V_1(t-s)P_K\tilde{F}_{-1}(s)\|_{L^q}
&
\le C\left\|V\left(\frac{t-s}{2}\right)e^{(\frac{t-s}{2})\Delta}P_K\tilde{F}_{-1}(s)\right\|_{L^q}
\\
&
\le C(t-s)^{-\frac{N}{2}(1-\frac{1}{q})}\|e^{(\frac{t-s}{2})\Delta}P_K\tilde{F}_{-1}(s)\|_{L^1}
\end{split}
\end{equation*}
for all $t\ge s+\delta>0$.
This implies that
\begin{equation}
\label{eq:6.19}
\begin{split}
t^{\frac{N}{2}(1-\frac{1}{q})}\tilde{J_4}(t)
&
\le t^{\frac{N}{2}(1-\frac{1}{q})}\int_0^{t/2}\|V_1(t-s)P_K\tilde{F}_{-1}(s)\|_{L^q}\,ds
\\
&
\le Ct^{\frac{N}{2}(1-\frac{1}{q})}
\int_0^{t/2}(t-s)^{-\frac{N}{2}(1-\frac{1}{q})}\|e^{(\frac{t-s}{2})\Delta}P_K\tilde{F}_{-1}(s)\|_{L^1}\,ds
\\
&
\le C\int_0^{t/2}\|e^{(\frac{t-s}{2})\Delta}P_K\tilde{F}_{-1}(s)\|_{L^1}\,ds,\qquad t>0.
\end{split}
\end{equation}
For any $T>0$, applying Lemma~\ref{Lemma:3.1}~(ii) with \eqref{eq:2.8},
we have
$$
\lim_{t\to\infty}(t-s)^{\frac{K}{2}}\|e^{(\frac{t-s}{2})\Delta}P_K\tilde{F}_{-1}(s)\|_{L^1}=0
$$ 
for any $s\in(0,T)$.
Then, by the Lebesgue dominated convergence theorem we see that
\begin{equation}
\label{eq:6.20}
\begin{split}
&
\limsup_{t\to\infty}\,t^{\frac{K}{2}}\int_0^T\|e^{(\frac{t-s}{2})\Delta}P_K\tilde{F}_{-1}(s)\|_{L^1}\,ds
\\
&
\le
\limsup_{t\to\infty}\,\int_0^T(t-s)^{\frac{K}{2}}\|e^{(\frac{t-s}{2})\Delta}P_K\tilde{F}_{-1}(s)\|_{L^1}\,ds
=0.
\end{split}
\end{equation}
Furthermore, 
applying Lemma~\ref{Lemma:3.1}~(i) with \eqref{eq:2.8},
for any $\delta>0$,
we deduce from \eqref{eq:3.7} that
\begin{equation}
\label{eq:6.21}
\begin{split}
\|e^{(\frac{t-s}{2})\Delta}P_K\tilde{F}_{-1}(s)\|_{L^1}
\le C(t-s)^{-\frac{K}{2}}\|P_K\tilde{F}_{-1}(s)\|_{L^1_K}
\le C(t-s)^{-\frac{K}{2}}E_{K,q}[\tilde{F}_{-1}](s)
\end{split}
\end{equation}
for all $t\ge s+\delta>0$.
Then, for any $T_0>0$, by \eqref{eq:6.16} and \eqref{eq:6.21} we see that
\begin{equation}
\label{eq:6.22}
\begin{split}
\int_T^{t/2}\|e^{(\frac{t-s}{2})\Delta}P_K\tilde{F}_{-1}(s)\|_{L^1}\,ds
&
\le C\int_T^{t/2}(t-s)^{-\frac{K}{2}}E_{K,q}[\tilde{F}_{-1}](s)\,ds
\\
&
\le Ct^{-\frac{K}{2}}\int_T^t(1+s)^{\frac{K}{2}-A-1}\,ds
\end{split}
\end{equation}
for all $t\ge2T$ and $T\ge T_0$.
Therefore, by \eqref{eq:6.19}, \eqref{eq:6.20} and \eqref{eq:6.22},
for any $d>0$ and $T\ge T_0$, we have
\begin{equation}
\label{eq:6.23}
t^{\frac{N}{2}(1-\frac{1}{q})}\tilde{J}_4(t)
\le dt^{-\frac{K}{2}}+C_*t^{-\frac{K}{2}}\int_T^t(1+s)^{\frac{K}{2}-A-1}\,ds
\end{equation}
for all sufficiently large $t$,
where $C_*$ is a constant independent of $T\in[T_0,\infty)$ and $d>0$.
Hence, by \eqref{eq:6.17}, \eqref{eq:6.18} and \eqref{eq:6.23} we obtain \eqref{eq:2.17} with $j=0$.

Next we assume that \eqref{eq:2.17} holds for some $j=m\in{\mathbb N}_0$.
Then, for any $q\in[1,\infty]$,
by \eqref{eq:1.3}, \eqref{eq:2.11}, \eqref{eq:2.16} with $j=m$ and \eqref{eq:6.14}
we have
\begin{equation}
\label{eq:6.24}
\begin{split}
(1+t)^{\frac{N}{2}(1-\frac{1}{q})}\|\tilde{F}_{m}(t)\|_{L^q}
&
\le C(1+t)^{\frac{N}{2}(1-\frac{1}{q})}\|u(t)-U_m(t)\|_{L^q}\times
\\
&\qquad\qquad
\times
\max\{\|u(t)\|_{L^\infty}^{p-1},\|U_m(t)\|_{L^\infty}^{p-1}\}
\\
&
\le C(1+t)^{\frac{N}{2}(1-\frac{1}{q})-A-1}\|u(t)-U_m(t)\|_{L^q},\qquad t>0
\end{split}
\end{equation}
Similarly to \eqref{eq:6.24} with \eqref{eq:2.10}, \eqref{eq:2.12} and \eqref{eq:2.13},
for any $\gamma\in\{1,\infty\}$,
we obtain
\begin{equation}
\label{eq:6.25}
\begin{split}
&
(1+t)^{\frac{N}{2}(1-\frac{1}{\gamma})}\|\tilde{F}_{m}(t)\|_{L^\gamma_K}
\\
&
\le C(1+t)^{\frac{N}{2}(1-\frac{1}{\gamma})}\|u(t)-U_m(t)\|_{L^\infty}
\times
\\
&
\qquad\qquad
\times
\max\{\|u(t)\|_{L^\infty}^{p-2}\|u(t)\|_{L^\gamma_K},
\|U_m(t)\|_{L^\infty}^{p-2}\|U_m(t)\|_{L^\gamma_K}\}
\\
&
\le C(1+t)^{\frac{N}{2}(1-\frac{1}{\gamma})+\frac{K}{2}-A-1}\|u(t)-U_m(t)\|_{L^\infty},
\qquad t>0.
\end{split}
\end{equation}
By \eqref{eq:2.17} with $j=m$, \eqref{eq:3.5}, \eqref{eq:6.24} and \eqref{eq:6.25} we see that
\begin{equation*}
\begin{split}
&
E_{K,q}[\tilde{F}_m](t)
\\
&
\le
C(1+t)^{\frac{K}{2}-A-1}\bigg((1+t)^{\frac{N}{2}(1-\frac{1}{q})}\|u(t)-U_m(t)\|_{L^q}
+\|u(t)-U_m(t)\|_{L^\infty}\bigg)
\\
&
=
\left\{
\begin{array}{ll}
o(t^{-A-1})+O(t^{\frac{K}{2}-(m+2)A-1})  & \mbox{if}\quad (m+1)A\not=K/2,\vspace{3pt}\\
O(t^{-A-1}\log t)  & \mbox{if}\quad (m+1)A=K/2,\\
\end{array}
\right.
\end{split}
\end{equation*}
as $t\to\infty$.
Therefore, applying same arguments as in the proof of the case $j=0$,
we can easily show that
$$
t^{\frac{N}{2}(1-\frac{1}{q})}\bigg(\tilde{J_2}(t)+\tilde{J_3}(t)+\tilde{J_4}(t)\bigg)
=
\left\{
\begin{array}{ll}
o(t^{-\frac{K}{2}})+O(t^{-(m+2)A})  & \mbox{if}\quad (m+2)A\not=K/2,\vspace{3pt}\\
O(t^{-\frac{K}{2}}\log t)  & \mbox{if}\quad (m+2)A=K/2,\\
\end{array}
\right.
$$
as $t\to\infty$.
This together with \eqref{eq:6.13} and \eqref{eq:6.15} implies \eqref{eq:2.17} with $j=m+1$.
Hence, by induction we see that \eqref{eq:2.17} holds for all $j\in{\mathbb N}_0$,
and the proof of Theorem~\ref{Theorem:2.2} is complete.
$\Box$
\vspace{10pt}

\noindent
{\large {\bf Acknowledgment.}}
The work of the first author (T. Kawakami) was supported in part 
by Grant-in-Aid for Young Scientists (B) (No. 24740107) and (No. 16K17629) 
of JSPS (Japan Society for the Promotion of Science)
and by the JSPS Program for Advancing Strategic International Networks 
to Accelerate the Circulation of Talented Researchers 
``Mathematical Science of Symmetry, Topology and Moduli, 
Evolution of International Research Network based on OCAMI'',
and the second author (H. Takeda) by Grant-in-Aid for for Young Scientists (B) (No. 15K17581) 
of JSPS.

\bibliographystyle{amsplain}

\end{document}